\documentclass{article}
\usepackage{latexsym, a4wide}
\usepackage{amsmath, rotating, color}
\usepackage{amsfonts,amssymb, amsthm}
\newcommand\unnumberedfootnote[1]{ %
        \let\temp=\thefootnote %
        \renewcommand{\thefootnote}{}%
        \footnote{#1}%
        \let\thefootnote=\temp%
        \addtocounter{footnote}{-1}}

\newtheorem{theorem}{Theorem}
\newtheorem{proposition}{Proposition}[section]
\newtheorem{lemma}[proposition]{Lemma}

\newtheorem{definition}[proposition]{Definition}
\theoremstyle{definition}
\newtheorem{remark}[proposition]{Remark}

\numberwithin{equation}{section}
\begin{document}
\title{\LARGE The process of most recent common ancestors\\  in an evolving coalescent}

\thispagestyle{empty}

\author{{\sc by P. Pfaffelhuber\thanks{Travel support from DFG,
      Bilateral Research Group
      FOR 498.} and A. Wakolbinger} \\[2ex]
  \emph{Ludwig-Maximilian University Munich} \\ \emph{and
    Goethe-University Frankfurt}\vspace*{-5ex}} \date{}

%\vspace*{-5ex}

\maketitle
\unnumberedfootnote{\emph{AMS 2000 subject classification.} 60K35
  (Primary) 92D25 (Secondary).}

\unnumberedfootnote{\emph{Keywords and phrases.} Kingman's Coalescent,
  Look-down process, most recent common ancestor}
\begin{abstract}
\noindent
Consider a haploid population which has evolved through an
exchangeable reproduction dynamics, and in which all individuals alive
at time $t$ have a most recent common ancestor (MRCA) who lived at
time $A_t$, say. As time goes on, not only the population but also its
genealogy evolves: some families will get lost from the population and
eventually a new MRCA will be established. For a time-stationary
situation and in the limit of infinite population size $N$ with time
measured in $N$ generations, i.e. in the scaling of population
genetics which leads to Fisher-Wright diffusions and Kingman's
coalescent, we study the process $\mathcal A = (A_t)$ whose jumps form
the point process of time pairs $(E,B)$ when new MRCAs are established
and when they lived. By representing these pairs as the entrance and
exit time of particles whose trajectories are embedded in the
look-down graph of Donnelly and Kurtz (1999) we can show by
exchangeability arguments that the times $E$ as well as the times $B$
from a Poisson process. Furthermore, the particle representation helps
to compute various features of the MRCA process, such as the
distribution of the coalescent at the instant when a new MRCA is
established, and the distribution of the number of MRCAs to come that
live in today's past. \end{abstract}

\section{Introduction}
%Consider a population of constant size whose individuals reproduce in an exchangeable way. 
The genealogy back to the most recent common ancestor (MRCA) of those
currently alive, and especially the time back to the MRCA, has been an
ongoing object of interest in mathematical population genetics, see
\cite{Littler1975}, \cite{Griffiths1980} for early references and
\cite {Wakeley2005} for a recent monograph. The limit of effective
population size $N\to \infty$, with time measured in units of $N$
generations, is the scaling in which Kingman's coalescent appears
(\cite{Kingman1982a}): in the rescaled time measured
backward from a fixed time $t$,  the number  of {\em
  ancestral lineages} enters from infinity and jumps
from $k$ to $k-1$ at rate $\binom{k}{2}$. (Here and below we assume
that the population size remains constant in time.) The depth $D_t$ of
the {\em coalescent tree}, that is the rescaled time it takes the number of ancestral lineages to decrease from $\infty$ to $1$, is then a sum of
exponentially distributed random variables  with mean
${\binom{k}{2}}^{-1}$, $k=2,3,\ldots$, and consequently has expectation
$2$.

With the population evolving further, also its genealogical relationships given by the coalescent tree change. In this study we are interested in
the time evolution of one particular characteristics of the genealogy,
that is, the time $A_t= t-D_t $ when the MRCA of the population at
time $t$ lived.  We will refer to $\mathcal A = (A_t)_{t\in\mathbb R}$
as the {\em MRCA process}.

At any time $t$ the total population consists of two oldest families,
which stem from the two oldest lines of descent dating back to the
MRCA who lived at time $A_t$. These two families will coexist for a
while after time point $t$, and during this time interval the path of
the MRCA process $\mathcal A$ stays constant. At some random time $E_t
> t$, one of the two families will go extinct and the other one will
fixate in the population. The MRCA of this surviving family must be
more recent than $A_t$, which amounts to a jump of the MRCA process at
time $E_t$.  Consequently, at time $E_t$, the next MRCA is
established, and the time when this next MRCA lives is $B_t :=
A_{E_t}$. In other words, the path of the process $\mathcal A$ is
constant as long as the two currently oldest families coexist in the
population, and jumps from $A_t$ to $B_t$ at time $E_t$ when one of
the two families fixates.

\medskip

The MRCA process is embedded in the genealogy of the population which
is assumed to evolve in a time stationary way. For a finite population
consisting of $N$ individuals, a way to construct the genealogy comes
with the graphical representation of the Moran model: for each ordered
pair $(i,j)$ of indices $i \neq j \in \{1,..,N\}$, an exponential
clock rings at rate $1/2$, and whenever this happens, the individual
with index $j$ dies and is replaced by an offspring of the individual
with index $i$.  This results in a partitioning of $ \mathbb R \times
\{1,\ldots, N\}$ into coalescing ancestral lineages, from which one
can read off a version of the MRCA process for $N$ individuals. This
process clearly inherits time stationarity from the evolution of the
population.

The Moran dynamics is exchangeable with respect to the individuals'
indices.  In contrast, the look-down process introduced by Donnelly
and Kurtz (1999), which is the basic tool in our study and will be
reviewed in Section \ref{2}, arranges the individuals' indices
(henceforth referred to as {\em levels}) at any time according to the
persistence of the individuals' offspring in the population: the
offspring of an individual at level $i$ outlives the offspring of any
contemporary individual at some higher level.  For a finite population
number $N$ this is achieved as follows: Each level $j$ ``looks down''
to each smaller level $i$ at rate 1. Whenever this happens, all
individuals at levels $j,\ldots, N-1$ are pushed one level up, the
individual at level $N$ is killed, and the individual at level $i$
spawns a child at level $j$. The time stationary MRCA process read off
from the look-down graph obviously has the same distribution as the
time stationary MRCA process read off from the Moran graph.

The look-down process allows a passage to the limit of infinite
population size in which the ordering by persistence is preserved. The
construction of the {\em random look-down graph} on $\mathbb R \times
\mathbb N$ proceeds in the very same way as described above, except
that there is no killing of individuals at any finite level. Instead,
the offspring of an individual at level $i\ge 2$ goes to extinction as
soon as this {\em line of ascent} of the individual is pushed to
infinity.  All this will be explained in more detail in Section
\ref{2}.

Because of the ordering by persistence, each MRCA of the population
lives at level 1 at some time $B$ at which it gives birth to an
individual at level 2. As soon as the offspring of these two
individuals fixates in the population, the MRCA is established. Again,
because of the ordering by persistence, this happens at the time $E$
when the line of ascent which was pushed at time $B$ from level 2 to
level 3 reaches infinity. The process $\mathcal F := \{(E,B)\}$, which
consists of all pairs of time points when an MRCA is established in
the population and when it lived, is a time-stationary point process;
we call it the {\em MRCA point process}. The paths of $\mathcal A$ and
the point configurations of $\mathcal F$ are in an obvious one-to-one
correspondence.
\begin{figure}
% \beginpicture
% \setcoordinatesystem units <.7cm,1cm>
% \setplotarea x from -.5 to 5, y from -1 to 4
% \plot 1 0 1 4 /
% \plot 5 0 5 4 /
% \plot 1 2.6 5 3.8 /

% \multiput {{\color{black}{\tiny $\bullet$}}} at 1 2.6 *100 0.04 .012 /

% %\setdots<1.4pt>
% %\plot 1 0.5 5 2.8 /

% \multiput {{\color{black}{\tiny $\bullet$}}} at 1 0.5 *20 0.2 .115 /

% \setdots
% \plot 1 2.8 5 2.8 /
% \put{$t$} [lC] at 5.2 2.8
% \put{$E_t$} [lC] at 5.2 3.8
% \put{next MRCA $B_t$} [rC] at 0.8 2.6
% \put{$\bullet$} [cC] at 5 3.8
% \put{$\bullet$} [cC] at 1 2.6
% \put{today's MRCA $A_t$} [rC] at 0.8 0.5
% \put{(a)} [cC] at 3 -0.3

% \setsolid
% \plot 11 0 11 4 /
% \plot 15 0 15 4 /

% \plot 11 0.2 15 1.2 /
% \multiput {{\color{black}{\tiny $\bullet$}}} at 11 0.2 *100 0.04 .01 /

% \put{$\bullet$} [cC] at 11 0.2
% \put{$\bullet$} [cC] at 15 1.2
% \put{$B$} [rC] at 10.8 0.2
% \put{$E$} [lC] at 15.2 1.2
% \plot 11 1.5 15 2.2 /
% \multiput {{\color{black}{\tiny $\bullet$}}} at 11 1.5 *100 0.04 .007 /

% \put{$\bullet$} [cC] at 11 1.5
% \put{$\bullet$} [cC] at 15 2.2
% \put{$B'$} [rC] at 10.8 1.5
% \put{$E'$} [lC] at 15.2 2.2
% \plot 11 2 15 3.8 /
% \multiput {{\color{black}{\tiny $\bullet$}}} at 11 2 *100 0.04 .018 /

% \put{$\bullet$} [cC] at 11 2
% \put{$\bullet$} [cC] at 15 3.8
% \put{$B''$} [rC] at 10.8 2
% \put{$E''$} [lC] at 15.2 3.8
% \put{(b)} [cC] at 13 -0.3

% \setdots
% \plot 15 1.2 11 1.2 /
% \plot 15 2.2 11 2.2 /
% \plot 15 3.8 11 3.8 /
% \endpicture
\begin{center}
\includegraphics[width=13.5cm]{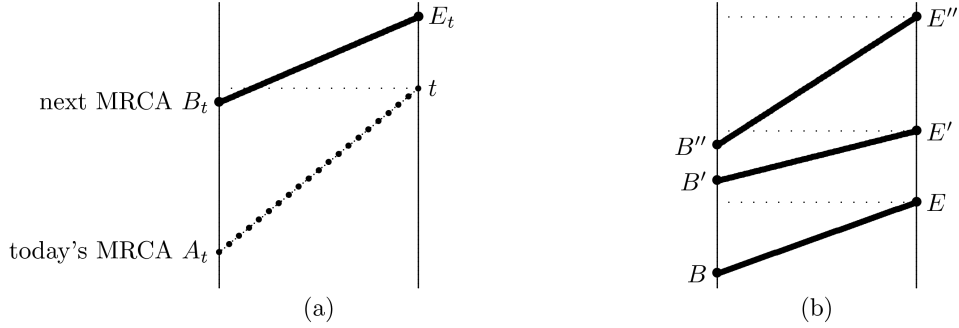}
\end{center}
\caption{\footnotesize\label{fig:PP} (a) At time $A_t$ the MRCA of the
  population at time $t$ lived, $E_t$ and $B_t$ are the times when the
  next MRCA is established and when it lived. \quad (b) MRCAs occur in
  a time-stationary manner. The dots on the $B$-axis are time points
  at which MRCAs lived. The dots on the $E$-axis are time points at
  which the MRCA changes.}
\end{figure}

The step from time $t$ to the next MRCA, which is established at time
$E_t$ and lives at time $B_t$, and an illustration of the MRCA point
process $\mathcal F$ are depicted in Figures \ref{fig:PP}(a) and
\ref{fig:PP}(b) respectively. In both Figures, the left axis contains
the times when MRCAs live, and the right axis gives the times when
MRCAs are established. The joint distribution of $E_t$ and $B_t$ will
be given in Theorem \ref{th1} in Section \ref{3}.  Figure
\ref{fig:PP}(b) displays part of the MRCA point process $\mathcal F$.
Remarkably, not only the points $B$ but also the points $E$ form a
time stationary Poisson process, see Theorem \ref{Th:full} in Section
\ref{4}. This will be proved by representing the times $(B,E), (B',
E'),\ldots$ as the entrance and exit times of particles: the trajectory
of a particle is attached to the line of ascent which is pushed from
level 2 to level 3 at time $B$ and exits at time $E$.  We will specify
the Markovian dynamics of this particle system, compute its
equilibrium distribution and show that, whenever a particle exits at
some time $E$, at this very instant the system of (remaining)
particles is in equilibrium. This allows to conclude that the waiting
time to the next exit time is exponential. The processes $\mathcal A$
and $\mathcal F$, however, are not Markov, see Remark
\ref{remarkBurke}.3.

In Theorem 3 we compute the distribution of the random
number $Z_t = \#\{(E,B) \in \mathcal F \, | \, E>t, B<t \}$ of MRCAs
that are established after time $t$ and live before time $t$. In
particular, it turns out that the probability that the next MRCA lives
in today's future is $\mathbf P[Z_t = 0] = 1/3$.

%In contrast to the MRCA particle system, the processes $\mathcal A$
%and $\mathcal F$ are not Markov, see Remark \ref{remarkBurke}.3.
%However, their past and future are independent given $\{Z_t = 0\}$.

As noted by \cite{Tajima1990}, the amount of polymorphism in a
population is related to the fixation of alleles. When an allele
fixates, the MRCA of the population must have changed. At such a
fixation time, the full coalescent is unusually short. As neutral
mutations fall independently on the branches of the genealogical tree,
this means that the amount of polymorphism is low at fixation times.

We start out by reviewing the look-down process in Section \ref{2},
describe our results in Sections \ref{3} and \ref{4}, point out some
relations to population genetics in Section \ref{6} and give the
proofs of Theorems 1-3 in Sections \ref{7}--\ref{9}.

\section{The MRCA process: a look-down construction}
\label{2}
At every time a continuum population which follows a Wright-Fisher (or
Fleming-Viot) dynamics has a genealogy given by Kingman's coalescent.
The look-down process introduced by Donnelly and Kurtz
(\cite{DonnellyKurtz1999}) not only gives a countable representation
of evolving allele frequencies but at the same time stores
genealogical relationships of all the individuals alive in the
population at all times. Consequently the MRCA process can be read off
from the look-down process.

\subsubsection*{The look-down graph: ancestral lineages, lines of
  ascent and ordering by persistence}
\begin{figure}
\label{lookdowngraph}
\begin{center}
\includegraphics[width=13.5cm]{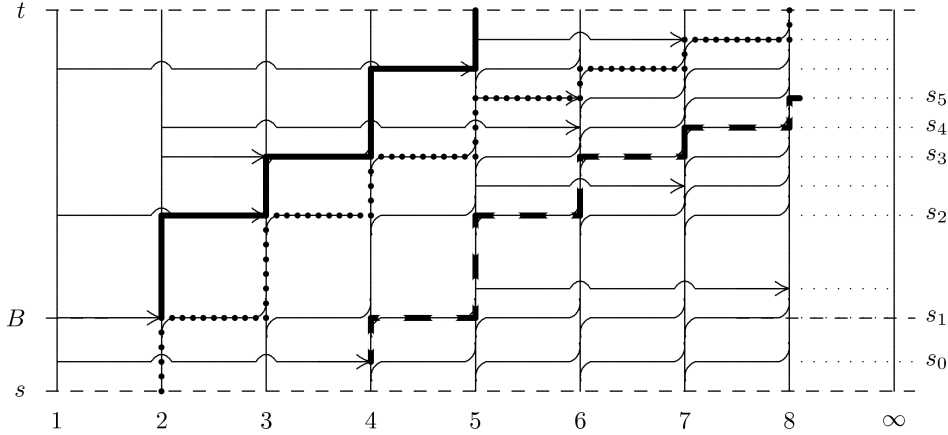}
\end{center}
\caption{\label{FigLookdown0}{\footnotesize Detail of a look-down
    graph. Time is running upwards; all lines at the first 8 levels
    are drawn between times $s$ and $t$. At times in $\mathcal P_{ij}$
    an arrow is drawn from $i$ to $j$. All lines at levels at and
    above $j$ are pushed upwards as indicated by bent lines. The solid
    marked line is the fixation curve $F_B^\tau, \tau \ge B$.  The
    dotted line is the coalescent curve $C_\tau^t(8), \tau \le t$. The
    dashed line is the line born at level 4 by the individual
    $(s_0,1)$; it is pushed one level up at times $s_1, s_2,\ldots$.
    In this picture, $X_s^t(1)=\ldots = X_s^t(5)=X_s^t(7)= 1$ and
    $X_s^t(6)= X_s^t(8)=2$; $Y_s^t(1)=1, Y_s^t(2)=6, Y_s^t(3)>8$;
    $C_s^t(1) = \ldots = C_s^t(5)=1$ and $C_s^t(6)=\ldots =
    C_s^t(8)=2$.}  }
\end{figure}

We first give a brief review of the ``modified look-down process''
(\cite{DonnellyKurtz1999}); see Figure \ref{FigLookdown0} for a
graphical illustration.

Consider the set of vertices
$$\mathcal V := \mathbb R \times \mathbb N.$$
We will refer to the vertex $(t,i)$ as the {\em individual at time}
$t$ {\em at level} $i$. For each ordered pair of levels $i<j$, let
$\mathcal P_{ij}$ be the support of a (rate one) Poisson point process
on $\mathbb R$, all these processes being independent. (In the
terminology of Donnelly and Kurtz, at each time $t \in \mathcal
P_{ij}$, the level $j$ looks down to level $i$.)  Based on the
processes $\mathcal P_{ij}$ we will construct a random countable
partition $\mathcal G$ of $\mathcal V$, whose partition elements will
be called \emph{lines}. The partition $\mathcal G$ will always contain
the so-called {\em immortal line} $\iota$ defined by
\begin{equation}\label{immline}
 \iota := \mathbb R \times \{1\}.
\end{equation}
For each $j > 1$, any point $s_0\in\bigcup_{i}\mathcal P_{ij}$ initiates a line
$G\in\mathcal G$ of the form
\begin{equation}\label{line}
  G = ([s_0,s_{1})\times\{j\}) \cup ([s_{1},s_{2})\times\{j+1\}) \cup ([s_{2},s_{3})\times\{j+2\})\cup \ldots
\end{equation}
with $s_{k+1}>s_k$ for all $k$. For a line $G$ as in \eqref{line} with
$s_0 \in \mathcal P_{ij}$ we say that $G$ is {\em born at level} $j$
{\em by the individual} $(s_0,i)$. We further say that $G$ is {\em
  pushed (one level up)} at times $s_{1},s_{2},\ldots$ and {\em
  exits} at time $s_{\infty}(G):= \lim_{n\to \infty} s_n$. The times
$s_k$ are given for $k=1,2,\ldots$ by $s_k = \inf\{s>s_{k-1}: s\in
\bigcup_{1\le \ell < m \le j + k-1} \mathcal P_{\ell m}\}$.
% Subsets of $ \mathcal V$ of the form \eqref{immline} or \eqref{line}
% will be called {\em lines}.  The random partition $\mathcal G$ is
% constructed from a countable family of Poisson processes:
Thus, a new line is born at level $j$ at each time $t$ when level $j$
looks down to some level $i < j$. Simultaneously, all the lines having
occupied at time $t-$ the levels $j, j+1,\ldots$ are pushed one level up.
Note that, since the pushing rate increases quadratically in $j$, the
exit time $s_\infty (G)$ is finite a.s.

For each $v \in \mathcal V$, we denote by $G_{v}$ the (unique) element
of $\mathcal G$ that contains $v$.  The {\em forward level process}
$Y_s^t(i),\, t\ge s$, initiated by the individual $u=(s,i)$ is given
by
$$Y_s^t(i):= \mbox{ level of } G_{u} \mbox{ at  time } t.$$
The {\em line of ascent} of individual $u=(s,i)$ is the part of line
$G_{u}$ after time $s$, that is
$$(t,Y_s^t(i))_{s\le t < s_\infty(g)}.$$

We say that a line $H$ {\em descends} from a line $G$ if either $H=G$,
or there is a finite sequence of lines $G_1,\ldots, G_{n-1} \in \mathcal
G$ such that $G_{k}$ is born by an individual in $G_{k-1}$ ,
$k=1,\ldots,n$, where $G_0:=G$ and $G_n := H$.

The {\em backward level process} $X_s^t(j),\, s\le t$, of an
individual $v=(t,j)$ arises by tracing back the level of $G_{v}$ to
the birth time of $G_{v}$, then jumping to the level of the individual
$u$ from which $G_{v}$ was born and tracing back the level of $G_{u}$
to the birth time of $(s,i)$, and so on.

The {\em ancestral lineage} of the individual $v=(t,j)$ is
$$ (s,X_s^t(j))_{-\infty < s \le  t}\; ;$$
note that eventually all ancestral lineages coalesce with the immortal
line.

We say that an individual $v \in \mathcal V $ {\em descends} from an
individual $u \in \mathcal V $ (or equivalently, $u$ is an {\em
  ancestor} of $v$) if $u$ belongs to the ancestral lineage of $v$.
The random tree spanning $\mathcal V$ which is obtained in this way is
the random {\em look-down graph}.
%; see Figure \ref{FigLookdown0}.

%\subsubsection*{Offspring trees and ordering by persistence}
%The {\em offspring tree} $\mathcal T_u$ of an individual $u$ consists
%of all individuals $v \in \mathcal V$ descending from $u$. The exit
%time $s_\infty (\mathcal T_u)$ of $\mathcal T_u$ is defined as the
%supremum of all times $t$ for which $\mathcal T_u$ contains an
%individual $(t,j)$ for some $j \in \mathbb N$.

Let $u=(s,i)$ and $v=(s,j)$ be two individuals living at the same time
$s$, with $i < j$. By construction the line of ascent of $v$ is pushed
whenever the line of ascent of $u$ is pushed, hence $Y_s^t(i) <
Y_s^t(j)$, and the line of ascent of $v$ exits not later than that of
$u$. In this sense, the {\em ordering of lines by contemporaneous
   levels} is an {\em ordering by persistence}.
Note also that for all times $s<t$ and all levels $i \in \mathbb N$:
$$ Y_s^t(i)=\inf \{j\in\mathbb N: X_s^t(j)=i\}.$$
Thus, the time when an individual's line of ascent reaches infinity
marks the time at which the individual's offspring goes extinct. 

%Let $u=(s,i)$ and $v=(s,j)$ be two individuals living at the same time
%$s$, with $i < j$. By construction the line of ascent of $v$ is pushed
%whenever the line of ascent of $u$ is pushed, hence $Y_s^t(i) <
%Y_s^t(j)$, and the line of ascent of $v$ exits not later than that of
%$u$. In this sense, the {\em ordering of lines by contemporaneous
%  levels} is an {\em ordering by persistence}.  In particular the
%consequence for each individual $u$ is
%$$s_\infty (\mathcal T_u) = s_\infty (G_u).$$
%Note also that for all times $s<t$ and all levels $i \in \mathbb N$:
%$$ Y_s^t(i)=\inf \{j\in\mathbb N: X_s^t(j)=i\}.$$
% The {\em offspring tree} $\mathcal T_u$ of an individual $u$
% consists of the line of ascent $\ell_u$ together with all lines in
% $\mathcal G$ descending from $\ell_u$. Since by definition the
% starting level of a newborn line is always larger than the level of
% its mother individual, it is clear that the time of extinction of
% $\mathcal T_u$ coincides with the exit time of the line of ascent
% $\ell_u$.

Let us note in passing that the ordering by persistence is a main
distinction between the version of the look-down process developed in
\cite{DonnellyKurtz1999} and its precursor introduced in
\cite{DonnellyKurtz1996}. In the latter, the order by persistence is
only stochastic, that is, lines of ascent of contemporaneous
individuals at lower levels are longer ``in probability''. In the
modified look-down process of \cite{DonnellyKurtz1999}, explained and
employed in the present paper, this property holds almost surely.

\subsubsection*{Coalescent curves and fixation curves} 
For $t \in \mathbb R$ and $i \in \mathbb N$ the {\em coalescent tree}
$\mathfrak C^t(i)$ consists of the ancestral lineages of the
individuals $(t,1),, (t,i)$, i.e.  $X_s^t(1), \ldots, X_s^t(i)$ for
$s\leq t$, whereas the {\em full coalescent tree} $\mathfrak C^t$ is
made up of the ancestral lineages of all individuals living at time
$t$.  All these lineages eventually coalesce with the immortal line.
Since any pair of ancestral lineages coalesces at rate 1, $\mathfrak
C^t(i)$ and $\mathfrak C^t$ are distributed like Kingman's (finite
respectively infinite) coalescent. The number of lineages remaining at
time $s < t$ can be expressed as
\begin{equation}\label{defCst}
C_s^t(i):=\max\{ X_s^t(1), ...,  X_s^t(i)\},   \qquad 
C_s^t:=C_s^t(\infty):=\sup_{i\in\mathbb N} C_s^t(i). 
\end{equation}
In words, $C_s^t(i)$ is the number of time $s$-ancestors of the time
$t$-individuals at levels $1,\ldots, i$, and $C_s^t$ is the number of
time $s$-ancestors of the whole population at time $t$.  For fixed
$t$, we call $(C_s^t)_{s\leq t}$ the \emph{coalescent curve} in the
look-down graph back from time $t$. It is distributed like the death
process in Kingman's coalescent entering from infinity.

The time when the MRCA of the total population at time $t$ lived is
$$A_t := \sup\{s : C_s^t = 1\}.$$
All individuals at time $t$ descend either from individual
$(A_t,1)$ or from individual $(A_t ,2)$. At time $A_t$ a line must be born at level $2$, which is equivalent to  $A_t \in \mathcal
P_{12}$.
Denote the next point in $\mathcal P_{12}$ after $A_t$ by $B_t$:
$$B_t:= \min\{s \in \mathcal P_{12} : s > A_t\}.$$
The offspring of the two individuals $(B_t,1),(B_t,2)$ evolves
towards fixation in the population by pushing the line of ascent of
the individual $(B_t,3)$ towards infinity. The time $E_t$ when this line
of ascent exits equals the time when the offspring of the individual
$(A_t,2)$ is expelled by the offspring of $\{(B_t,1),(B_t,2)\}$. Thus the
time $E_t$ is the first time after $t$ when a new MRCA is established,
and the time when this MRCA lives is $B_t$.

Note that at any time $\tau$ between $B_t$ and $E_t$, all the levels
$1,\ldots, Y_B^\tau(3)-1$ are occupied by offspring of
$\{(B_t,1),(B_t,2)\}$, whereas level $Y_B^\tau(3)$ is not.  We
therefore call
\begin{align}\label{YBt} F_{B_t}^\tau:= Y_{B_t}^\tau(3)-1 =
  Y_{A_t}^\tau(2)-1, \quad B_t\le \tau < E_t,
\end{align}
the  {\em fixation curve} starting in time $B_t$. (For the
equality in \eqref{YBt}, note that the line containing $(B_t,3)$ was
born at time $A_t$ at level 2 and was pushed to level 3 at time $B$.)
When $Y_{B_t}^\tau(3)=k$ the corresponding line moves to $k+1$ at the next
look-down event among the first $k$ levels, i.e. at rate
$\binom{k}{2}$. As a consequence, $F_{B_t}^\tau$ is pushed from level $k$
to level $k+1$ at rate $\binom{k+1}{2}$.
  
The MRCA point process $\mathcal F$ records all the time points when
the fixation curves start and end. We will pursue this in Section
\ref{4}, by constructing an autonomous particle system whose
trajectories give the fixation curves.
  
Whereas the coalescent curves are constructed from any $t$ backwards
in time, the fixation curves start only at points in $\mathcal P_{12}$
and are constructed forwards in time.

At a time $E$ when a fixation curve ends (and a new MRCA is
established), all individuals descend from the MRCA who lived at the
time $B$ when this fixation curve started. Hence the fixation curve
between time points $B$ and $E$ equals the coalescent curve back from
time $E$.  With time proceeding, the coalescent curve evolves, being
more and more ``zipped away'' from the upper end of the fixation curve
(near time point $E$), and still sharing the lower part (near time
point $B$) for a while.

% Using either fixation curves or coalescent curves there are two
% different ways to describe the MRCA process.  Recall that $A_t$, the
% time when the MRCA of the individuals $(t,1), (t,2),..$ lived, is
% expressed in terms of the coalescent curve back from time $t$:
%  $$ A_t := \sup\{s\leq t: C_s =1\}$$
%  This equals the starting time of the most recent fixation curve which was completed not later than at time $t$:
%  \begin{equation}\label{twoways}
%A_t =\sup\{s\leq t: Y_s^t(2) \geq k\;  \forall k \in \mathbb N\}.
%\end{equation}
%In fact, for all $j=1,2,..$  one observes the identity
%\begin{equation}\label{twowaysj}
% A_t^j := \sup\{s\leq t: C_s^t(j)=1\}= \sup\{s\leq t: Y_s^t(2)-1\geq j\} =: \tilde A_t^j,
%\end{equation}
%from which \eqref{twoways} follows by taking $j \to \infty$.

Having now constructed the process $\mathcal A$ in terms of the
look-down graph, we will study its properties in the next sections.

\section{From today to the next MRCA}
\label{3}
As in the previous section, $A_t$ denotes the time when the current
MRCA lived, $E_t$ is the time when the next MRCA is established and
$B_t=A_{E_t}$ is the time when the next MRCA lives. In this section we
will compute the conditional distribution of $(E_t, B_t)$ given $A_t$.

The following random variables will play a crucial role:
\begin{equation}\label{defL}
L_t:= F_{B_t}^t\, ,
\end{equation}
the level at time $t$ of the fixation curve starting at time $B_t$,
and
\begin{equation}\label{defI}
I_t:= C_{B_t}^t\, ,
\end{equation}
the level at time $B_t$ of the coalescent curve back from time $t$,
where we define $$L_t:= 1 \mbox{ and } I_t:= \infty \mbox{ on the
  event }\{B_t > t\}.$$

Without loss of generality, and to ease notation, let us put $t=0$,
and write $L:= L_0, I:= I_0, E:= E_0, B:= B_0$.

Note that, because of the ordering by persistence, the lines of ascent
starting at time $0$ from levels $1,\ldots, L$ exit only after time $E$,
whereas the lines starting at time $0$ from levels $L+1, L+2,\ldots$ exit
at time $E$ or earlier. Thus, $L$ is the random number of individuals
in the present population that still have offspring when the next MRCA
is established.

\begin{proposition} \label{LI_prop} The pair $(L,I)$ is independent of
  $A_0$ and has distribution
  \begin{equation}\label{IL}
    \begin{aligned}
  \mathbf P[L=\ell, I=i] = \begin{cases}\displaystyle \frac{\ell-1}{3\binom{\ell+i}{\ell}},&\ell\geq 2, i\geq 3\\
    \tfrac 13, & \ell=1,i=\infty\\ 0, & \text{else.}\end{cases}
\end{aligned}\end{equation}
\end{proposition}
Proposition \ref{LI_prop} will be proved in Section \ref{7}. 
\begin{remark} \label{LIremark}
\begin{enumerate}
\item Summing over $i$ in \eqref{IL} leads to the distribution of
    $L$:
    \begin{align}\label{distL1}
    \mathbf P[L=\ell]  =
    \frac{2}{(\ell+1)(\ell+2)}, \quad \ell=1,2,\ldots
  \end{align}
  Since $\{L=1\} = \{B>0\}$ is the event that that first fixation
  curve which ends after time $0$ has not yet started by time $t=0$.
  we infer that the probability that the next MRCA lives in today's
  future is
  $$\mathbf P [B>0] = \mathbf P [L=1] = 1/3.$$
\item Here is another quick way to \eqref{distL1}, exploiting
  exchangeability.  Recall that the number $L$ gives the number of
  lines that still have offspring at the time when the next MRCA is
  established. At any time there are two oldest families in the
  population. The family sizes of these two oldest families, denoted
  by $P$ and $1-P$, evolve according to a Wright-Fisher diffusion.  It
  is well known (and can be understood from the P\'olya urn scheme
  embedded in the genealogy; see e.g. facts about the
  P\'olya-Eggenberger distribution in \cite{JohnsonKotz1977}, eq.
  (4.1)) that, at any fixed time, say at time $t=0$, $P$ is uniformly
  distributed on $[0,1]$.  This also remains true conditioned on the
  event $A_0^\infty=-d$. By exchangeability, the probability that the
  first $\ell$ most persistent lines are in one and the $(\ell+1)$-st
  most persistent line is in the other family is
  \begin{align}\label{distL}
    \mathbf P[L=\ell] = 2\int_0^1 p^\ell(1-p)dp = 2\Big(
    \frac{1}{\ell+1} - \frac{1}{\ell+2}\Big) =
    \frac{2}{(\ell+1)(\ell+2)}.
  \end{align}
  \end{enumerate}
\end{remark}

To prepare for Theorem \ref{th1}, we need one more bit of notation.
\begin{definition}
\label{defTk}
Let $T_k$ be independent exponentially distributed random variables
with parameter $\binom{k}{2}$, $k=2,3,\ldots$, and
\begin{equation*}
S_i^j = \sum_{k=i+1}^{j} T_k, \qquad 1\le i < j \le \infty.
\end{equation*}

For $d > 0$ and $i=1,2,..$ let $R_{i,d}$ be a random variable whose
distribution equals the conditional distribution of $S_i^\infty$ given
that $S_1^i + S_{i}^\infty =d.$
\end{definition}

The random variable $S_i^j$ represents the time which Kingman's
coalescent requires to come down from $j$ to $i$. Consequently,
$R_{i,d}$ refers to the random time for a coalescent to come down from
infinity to $i$ lines, given that coming down to 1 line requires
exactly time $d$. Note also that $S_i^j$ represents the time a
fixation curve needs to be pushed from level $i$ to level $j$. This
can be seen because a fixation curve goes from level $\ell$ to level
$\ell+1$ whenever a look-down event among the first $\ell+1$ levels
occurs, i.e. with rate $\binom{\ell+1}{2}$.

We are now prepared to state Theorem \ref{th1}, which together with
Proposition \ref{LI_prop} yields the desired conditional distribution
of $(E,B)$ given $A_0$.

\begin{theorem}
\label{th1}
Let $L$ and $I$ be as in \eqref{defL} and \eqref{defI}.  The
conditional distribution of $(E,B)$, given $A_0 = -d$, $L=\ell$ and $I
= i$ is represented by the random variables
\begin{equation}\label{dynamics}
\begin{aligned}
\begin{cases}  (S_{1}^2 + S_2^\infty, S_1^2) & \mbox{ if } \ell =1,\\[2ex]
 (S_{\ell}^\infty, -R_{i,d})
 &  \mbox{ if } \ell >1,\end{cases}
\end{aligned}\end{equation}
where $(S_{\ell}^j)$ and $R_{i,d}$ have the distribution specified in
Definition \ref{defTk}, and $R_{i,d}$ and $S_\ell^\infty$ are
independent.
\end{theorem}
The proof of Theorem \ref{th1} is given in Section \ref{7}.

\begin{remark}
  \label{remExp}
  \begin{enumerate}
  \item Combining \eqref{dynamics} and \eqref{distL1} we obtain
 \begin{equation}\label{expid}
    \mathbf P[E\in ds \big | A_0=d] = \sum_{\ell=1}^\infty \frac{2}{(\ell+1)(\ell+2)} 
    \mathbf P[S_\ell^\infty\in ds], \qquad s \ge 0.
  \end{equation}
  From this one can conclude that the conditional distribution of $E$
  given $A_0$ is standard exponential. Indeed, think of a 2-sample
  (i.e. a subsample of size two) embedded in a full coalescent.  The
  coalescence time of this 2-sample is standard exponentially
  distributed.  Denoting by $L'$ the number of lineages remaining in
  the full coalescent at the time when the 2-sample has found its
  common ancestor, one sees from \cite{GriffithsTavare2003}, eq.
  (2.10) or \cite{Saundersetal1984}, Lemma 3, or by direct
  calculation, that $L'$ has the same distribution as $L$ specified in
  \eqref{distL1} and \eqref{distL}.  This shows that the r.h.s. of
  \eqref{expid} is a decomposition of the standard exponential
  distribution.
\item Here is another quick (though slightly informal) argument that
  the waiting time to the next jump of the MRCA is exponential,
  independently of the depth of the current MRCA. Note first that,
  conditioned on $A_0=-d$ the split of the population size into the
  two oldest families at time $t=0$ is uniformly distributed on
  $[0,1]$. As a consequence, given the MRCA does not jump during the
  time interval $[0,s]$, the split remains uniformly distributed also
  at time $s$. (This corresponds to the fact that the uniform
  distribution is a quasi-equilibrium for the Wright-Fisher
  diffusion.)  At the next jump of the MRCA process one of the two
  oldest families dies out.  After the jump there will be two families
  inside the surviving family that again make up a uniform split.
  This implies that the time between jumps proceeds in a memoryless
  manner, showing that the conditional distribution of $E_0$ given
  $A_0$ is exponential.

  Notably, the fact of exponential waiting times between the jumps can
  also be read from (3.10) in \cite{Watterson1982a}. See Section
  \ref{6} for comments relating to this paper and to other
  applications.
 \end{enumerate}
\end{remark}

\section{A particle representation of the MRCA point process}
\label{4}
The set $\mathcal G$ of lines defined in Section \ref{2} randomly
partitions the set $\mathcal V = \mathbb R \times \mathbb N$. Let us
write
\begin{equation*}
  \mathcal G_2 := \{G \in \mathcal G \, \big | \, G \mbox{ is born at level } 2 \}.
\end{equation*}
For each line $G \in \mathcal G_2$ we write $B:= B(G)$ for the time
when $G \in \mathcal G_2$ is pushed from level 2 to level 3 (due to
the birth of the next line in $\mathcal G_2$) and $E:= E(G)$ for the
exit time of $G$. Thus we obtain a one-to-one correspondence between $
\mathcal G_2$ and the sequence of fixation curves by associating with
any $G \in \mathcal G_2$ the fixation curve $F_B$ starting at time $B$
and ending at time $E$. This fixation curve is related to the level
path of $G$ by $F_B^\tau = Y_B^\tau(3)-1$ for $B \le \tau < E$; see
\eqref{YBt}. The MRCA point process $\mathcal F$ then can be written
as
\begin{equation*}
  \mathcal F = \{(E,B)\, \big | \, G \in \mathcal G_2\}.
\end{equation*}
Additionally, we write
\begin{equation*}
  \eta := \{E\, | \, (E,B) \in \mathcal F\}\qquad\text{ and }\qquad \eta_t := \eta \cap (-\infty, t]
\end{equation*} 
for the exit time point process and its restriction to $(-\infty, t]$
respectively.

In this section we will gain more information about the processes
$\mathcal F$ and $\eta$ by interpreting the fixation curves as the
trajectories of an interacting particle system on $\{2,3,4,\ldots\}$
whose dynamics and equilibrium distribution we will compute.

Let 
\begin{equation*}
  Z_t:= \#\{(E,B) \in \mathcal F \, \big | \, E>t, B<t \}.
\end{equation*}
In other words, $Z_t$ is the number of fixation curves present at time
$t$, that is, the number of MRCAs which will be established after time
$t$ and have lived before time $t$.

Write
  \begin{equation}\label{particles}
L_t^1 > L_t^2 > \ldots > L_t^{Z_t}>1
\end{equation}
for the levels of the fixation curves at time $t$. Let us interpret
$(L_t^1, L_t^2, \ldots, L_t^{Z_t})$ as a configuration of particles on
the set of levels $\{2,3,4,\ldots\}$ at time $t$, and put
$$\Lambda_t := (L_t^1, L_t^2,\ldots ),\qquad \Lambda := (\Lambda_t),$$
where $L_t^j := 1$ for $j > Z_t$. The first components in the MRCA
point process $\mathcal F$ are the exit times of the ``leading
particles'', i.e. those time points $E$ where $\lim_{t\uparrow E}L_t^1
= \infty$. Whenever a particle exits, the indices of all remaining
particles are shifted down by one:
\begin{equation}\label{jumpback}
(L_E^1, L_E^2,\ldots ):= (L_{E-}^2, L_{E-}^3,\ldots ).
\end{equation}
Here is a verbal description of the dynamics of the particle system
(see Proposition \ref{propLambda} for a formal statement): Particles
are pushed in at level 2 at rate 1, each particle at level $\ell \ge
2$ is pushed one level up at rate $\binom{\ell + 1}{2}$, and this is
done in a coupled way such that, whenever a particle is pushed, all
particles at higher levels are pushed simultaneously.  The next
theorem specifies the equilibrium distribution of $\Lambda$. We will
see that this distribution prevails also in the distinguished random
time points $E$ where $\lim_{t\uparrow E}L_t^1 = \infty$. This
property is crucial to see that $\eta$ is a Poisson process.
 
\begin{figure}
\begin{center}
\includegraphics[width=13.5cm]{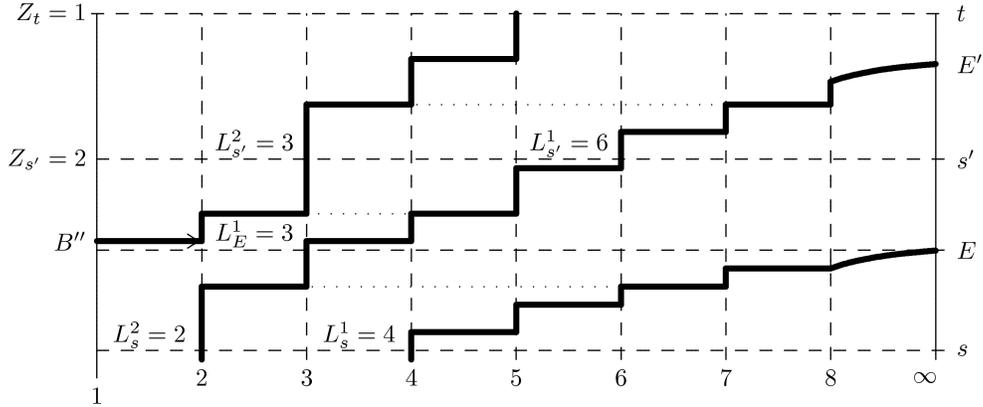}
\end{center}
\caption{\label{FigLookdown}The embedding of the fixation curves in
  the look-down process. At times $B$ fixation curves start and at
  times $E$ they end. In this example, at time $s$ the  number of particles in the system is $Z_{s'}=2$, the leading particle being at level $6$, and the second particle at level $2$.}
\end{figure}

%The equilibrium distribution of $\mathcal L$ can now be described.

\begin{theorem}
\label{Th:full} 
\begin{enumerate}
\item The process $\Lambda =(\Lambda_t)$ is Markov with stationary
  distribution
\begin{equation}\label{equlambda} \pi_{\Lambda} (\ell_1,\ell_2,\ldots) = \begin{cases}\displaystyle
    \frac 13 \prod_{j:\ell_j>1} \frac{2}{(\ell_j+2)(\ell_j-1)},&
    \ell_k>\ell_{k+1}, \quad \text{ if } \ell_{k}>1,\\ 0, &
    \text{else}.
  \end{cases}
  \end{equation} In particular, the stationary distribution of $L^1$ is
\begin{equation} \label{equilL}     \pi_{L^1}(\ell) = \frac{2}{(\ell+1)(\ell+2)}.
\end{equation}

%  In equilibrium, $L_t$ and $A_t$ are independent. At time points $E$
%  where $\lim_{t\uparrow E}L_t^1 = \infty$, the distribution of $L_E$
%  equals $\pi_{\Lambda}$, and $L_E$ and $A_E$ are independent.

%\item For all $t \in \mathbb R$, conditional under $C^t$, the particle
%configuration $\Lambda _t$ has distribution $\pi_\Lambda$.
%\item For all $t \in \mathbb R$, conditional under the coalescence
%  curve $(C^t_s)_{s\leq t}$, the past exit time points $\eta_t$ and
%  the event $\{t \in \eta\}$, the particle configuration $\Lambda _t$
%  has distribution $\pi_\Lambda$.
\item The process of exit times $\eta$ is a stationary Poisson
  process.
\end{enumerate}
\end{theorem}

\begin{remark}
\label{remarkBurke}
\begin{enumerate}
  % \item Yet another proof of the exponentiality of the waiting times
  %   between jumps of MRCAs (see Remark \ref{remExp}) is the
  %   following: At time 0 the process $\Lambda$ is in equilibrium,
  %   and so it is at time $t > 0$. Let $E$ be the first exit time of
  %   the leading particle after time 0. By Theorem \ref{Th:full}, the
  %   process $\Lambda$ is in equilibrium at time $E$, and hence also
  %   at time $t$, given the event $t>E$, and consequently also at
  %   time $t$, given the event $t<E$, which implies that the time $E$
  %   has an exponential distribution.

\item The ``arrival time points'' $B$ of the particles in the system $
  \Lambda$ (the times when the MRCAs live) are the points of the
  stationary Poisson process $\mathcal P_{12}$. Theorem \ref{Th:full}
  states that also the ``departure time points'' $E$ (the times when
  the MRCAs are established) form a stationary Poisson process.  Thus,
  the Poisson input process of times $B$ when the fixation curves
  begin is transformed by a ``dependent stochastic shift'' into the
  Poisson output process of times $E$ when they end. This is similar
  to Burke's theorem which states that the departure process in a time
  stationary $M/M/1$ queue is Poisson, see \cite{Kurtz1998} and
  references given there. A crucial property (proved already by Burke
  (1956)) is that in a stationary $M/M/1$ queue the distribution of
  the queue length at time $t$ is independent of the departure times
  $\le t$.  In the language of queueing theory, our particles
  correspond to customers entering at the time points of a stationary
  Poisson process, and the time which a typical customer spends in the
  system is distributed like $S_2^\infty$ specified in Definition
  \ref{defTk}.  These times are mutually dependent. As the proof of
  Theorem \ref {Th:full} reveals, like in Burke's theorem the state of
  the system (now given by the configuration $\Lambda_t$ of particles
  at time $t$) does not depend on the departure times $\le t$.

\item We recently learned from Tom Kurtz about the manuscript  \cite{DonnellyKurtz2006} where he and Peter Donnelly have
  established  the filtered martingale
  problem for the $N$-level analogue of the pair $(\Lambda, \eta)$ in
  the context of \cite{Kurtz1998}, Theorem 3.2 and thus achieved an
  alternative proof of the fact that the ``MRCA fixation process''
  $\eta $ is Poisson.

% \item Theorem \ref{Th:full} shows that not only the ``entrance time
%   points'' $B$ (the times when the MRCAs live) but also the ``exit
%   time points'' $E$ (the times when the MRCAs are established) from a
%   stationary Poisson process.  The Poisson input process of times $B$
%   when the fixation curves begin is transformed by a ``dependent
%   stochastic shift'' into the Poisson output process of times $E$ when
%   they end: new particles enter a system at Poissonian times, move in
%   a dependent way, and leave at Poissonian times.

% \item In Theorem \ref{Th:full} we are faced with a Markov process
%   which is reset to its equilibrium distribution along an increasing
%   sequence of stopping times.  This is closely related to Burke's
%   output theorem, see \cite{Kurtz1998} and references given there. We
%   recently learned from Tom Kurtz that he and Peter Donnelly have
%   established in \cite{DonnellyKurtz2006} the filtered martingale
%   problem for the $N$-level analogue of the pair $(\Lambda, \eta)$ in
%   the context of \cite{Kurtz1998}, Theorem 3.2, thus giving an
%   alternative proof of the fact that the ``MRCA fixation process''
%   $\eta$ is Poisson.

\item Whereas the particle process $\Lambda$ is Markov, the MRCA
  process $\mathcal A$ is not. This can be seen as follows:

\begin{figure}\label{Nonmarkov}
% \beginpicture
% \setcoordinatesystem units <.7cm,1cm>
% \setplotarea x from -7 to 5, y from -1 to 4
% \plot 1 0 1 4 /
% \plot 5 0 5 4 /

% \multiput {{\color{black}{\tiny $\bullet$}}} at 1 0.8 *100 0.04 .030 /
% \multiput {{\color{black}{\tiny $\bullet$}}} at 1 0.5 *20 0.2 .03 /
% \multiput {{\color{black}{\tiny $\bullet$}}} at 1 0.5 *20 0.2 .115 /

% \setdots
% \put{$t$} [lC] at 5.2 2.8
% \put{$s$} [lC] at 5.2 1.1
% \put{$E_t=E_s$} [lC] at 5.2 3.8
% \put{$B_t=B_s$} [rC] at 0.8 0.9
% \put{$\bullet$} [cC] at 5 3.8
% \put{$\bullet$} [cC] at 1 0.8
% \put{$A_t=A_s=a$} [rC] at 0.8 0.5
% \setdots
% \plot 1 0.5 7 0.5 /
% \plot 1 1.1 7 1.1 /
% \plot 1 2.8 5 2.8 /
% \setsolid
% \arrow  <0.2cm> [0.375,1] from 7 0.5 to 7 1.1
% \arrow  <0.2cm> [0.375,1] from 7 1.08 to 7 0.5
% %\arrow  <0.2cm> [0.375,1] from -1 1.12 to -1 2.78
% %\arrow  <0.2cm> [0.375,1] from -1 2.78 to -1 1.12
% %\arrow  <0.2cm> [0.375,1] from -1 3.8 to -1 2.82
% \put{$\varepsilon$} [lC] at 7.2 .85
% %\put{$\beta$} [rC] at -1.2 1.85
% %\put{$\gamma$} [rC] at -1.2 3.3
% \endpicture
\begin{center}
\includegraphics[width=7.5cm]{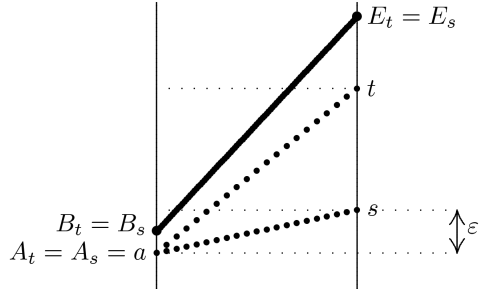}
\end{center}
\caption{\footnotesize\label{figNotMarkov}Assume we know $A_t=A_s=a$
  for some $s=a+\varepsilon$. This knowledge leads to a higher
  chance of the MRCA time $B_t$ falling between times $a$ and $s$ than
  in an equilibrium situation. This shows that the future of the process $\mathcal
  A$ at time $t$ depends on the past and $\mathcal A$ cannot be
  Markov. See text for explanation.}
\end{figure}

Let $a, s, t$ be as in Figure \ref{Nonmarkov}.  Conditioned on $A_t=a$
we obtain from Theorem 1:
$$ \mathbf P[B_t<s|A_t=a] = \mathbf P[R_{I,t-a}>t-s]
\xrightarrow{s\downarrow a}0.$$ On the other hand, we claim that $
\mathbf P[B_t<s|A_t=A_s=a]$ does not converge to zero as $s\downarrow
a$, which shows that $\mathcal A$ cannot be Markov. To verify the
claim, we write, using Bayes' rule
$$ \mathbf P[B_t<s|A_t=A_s=a] =
\frac{\mathbf P[B_s<s, A_t=A_s|A_s=a]}{\mathbf P[A_t=A_s|A_s=a]}.$$ By
Theorem \ref{th1}, the denominator converges to $e^{-(t-a)}$ as
$s\downarrow a$.  Likewise, the numerator is bounded away from $0$ as
$s\downarrow a$, a trivial lower bound being $$\mathbf P[L_s = 2]\,
\mathbf P[S_2^3 \ge t-s] = \frac 16 e^{-3(t-s)} \ge \frac 16
e^{-3(t-a)}.$$
\item The level $L_t^1$ of the leading particle in \eqref{particles}
  coincides with $L_t$ defined in \eqref{defL}. Thus we recover
  \eqref{distL1} from \eqref{equilL}.
\end{enumerate}
\end{remark}

\noindent 
Recall from \eqref{particles} that
\begin{equation*}
  Z_t = \max\{j \in \mathbb N \, \big | \, L_t^j > 1\},
\end{equation*}
where $\max \emptyset := 0$. Consequently, $$\{Z_t = 0\} = \{L_t^1 =
1\}.$$ This is the event that there is no particle on $\{2,3,4,..\}$
at time $t$, or equivalently, that all fixation curves starting before
time $t$ also end before $t$.  Given this event, the fixation curves
starting before time $t$ are independent of those starting after $t$.

In a way, the random variable $Z_t$ of MRCAs that are established in
today's future and live in today's past measures the dependence
between past and future in the MRCA process. Note also that because of
Theorem 3 the distribution of $Z_t$ does not change when $t$ is
conditioned to be the time of an MRCA change.  

In the next theorem we calculate the equilibrium distribution of
$Z:=Z_t$.
\begin{theorem}
\label{LawZ}
\begin{enumerate}
\item The probability generating function of $Z$ is
$$ \mathbf E[u^Z] = \frac 13 \exp\Big( \sum_{i=2}^\infty \log\Big( \frac{i(i+1)+2(u-1)}{(i+2)(i-1)}\Big).$$
\item The expectation and variance are 
$$\mathbf E[Z] = 1,\qquad \mathbf{Var}[Z] = 14 - \tfrac 43 \pi^2 \approx 0.84052.$$
\item The probability weights are given by
$$ \mathbf P[Z=z] =
\frac{2^z}{3}\sum_{\underline a: \sum i a_i=z} (-1)^{z + \sum a_i} \frac{1}{a_1!\cdots a_z!} \prod_{j=1}^z \Big(\frac{x_j}{j}\Big)^{a_j},\qquad (z \geq 0)$$
where the $a_i \in \mathbb N_0$ and
\begin{equation}
\label{Def:xk}
\begin{aligned}
 x_k &:= \frac{(-1)^{k+1}}{3^{2k-1}} \sum_{j=1}^{k}\binom{2k-j-1}{k-j}3^{j-1}\big(b_j - 1_{\{j\; \rm{ even}\}}2\zeta(j)\big),\\
 \zeta(k) &:= \sum_{j=1}^\infty \frac{1}{j^k},\\
b_j&:= 1 + \frac{1}{2^j} + \frac{1}{3^j}.
\end{aligned}
\end{equation}
\item The weights for $z=0,1,2,3$ are
%\begin{align*}
%x_0&=1 &&& x_1 & =\frac{11}{18},\\
%x_2&=-\frac{31}{108}+\frac{1}{27}\pi^2,&&& x_3 &= \frac{809}{5832}-\frac{1}{81}\pi^2,\\
%x_4&= -\frac{7121}{104976} + \frac{10}{2187}\pi^2 + \frac{1}{3645}\pi^4
%&&& x_5 &=\frac{7009}{209952} - \frac{35}{19683}\pi^2 - \frac{1}{6561}\pi^4
%\end{align*}
\begin{align*}
  \mathbf P[Z=0] & = \frac 13, &&&
  \mathbf P[Z=1] & = \frac{11}{27}\approx 0.40740,\\
  \mathbf P[Z=2] & = \frac{107}{243}-\frac{2}{81}\pi^2\approx 0.19664,
  &&& \mathbf P[Z=3] & = \frac{1003}{2187}-\frac{10}{243}\pi^2\approx
  0.05246.
\end{align*}
%\vspace{-5ex}

%\begin{align*}
%  \mathbf P[Z=4] & = \frac{71125}{157464} - \frac{311}{6561}\pi^2 + \frac{8}{32805}\pi^4\approx 0.00761,\\
%  \mathbf P[Z=5] & = \frac{28201}{59049} - \frac{3526}{59049}\pi^2 +
%  \frac{38}{32805}\pi^4\approx 0.00108.
%\end{align*}

% \smallskip
% \noindent
% \begin{tabular}{|c|c|c|c|c|c|}\hline
% \raisebox{0ex}[3ex][2ex]{}$\ell$              & 0     & 1 & 2 & 3 & 4  \\\hline
% \raisebox{0ex}[3ex][2ex]{$\mathbf P[Z=\ell]=$} & $\frac{1}{3}$ & $\frac{11}{27}$ & $\frac{107}{243}-\frac{2}{81}\pi^2$ & $\frac{1003}{2187}-\frac{10}{243}\pi^2$ & $\frac{71125}{157464} - \frac{311}{6561}\pi^2 + \frac{8}{32805}\pi^4$  \\\hline
% \raisebox{0ex}[3ex][2ex]{}$\mathbf P[Z=\ell]\approx$ & 0.33333 & 0.40740 & 0.19664 & 0.05246 & 0.00761  \\\hline
% \end{tabular}
% \smallskip
%
% \noindent
% \begin{tabular}{|c|c|c|}\hline
% \raisebox{0ex}[3ex][2ex]{}$\ell$      & 5 & 6  \\\hline
% \raisebox{0ex}[3ex][2ex]{$\mathbf P[Z=\ell]=$} &  $\frac{28201}{59049} - \frac{3526}{59049}\pi^2 + \frac{38}{32805}\pi^4$ & $\frac{256625}{531441} - \frac{1282}{19683}\pi^2 + \frac{502}{295245}\pi^4 - \frac{4}{688905}\pi^6$\\\hline
% \raisebox{0ex}[3ex][2ex]{}$\mathbf P[Z=\ell]\approx$ & 0.00108 & 0.00010\\\hline
% \end{tabular}
%
% \noindent
% \begin{align*}
%  \mathbf P[Z=0] &= \frac{1}{3}\approx 0.33333,\\
%  \mathbf P[Z=1] &= \frac{11}{27}\approx 0.40701,\\
%  \mathbf P[Z=2] &= \frac{4}{162} \Big( \frac{107}{6} - \pi^2\Big)\approx 0.19664,\\
%  \mathbf P[Z=3] &= \frac{8}{972} \Big( \frac{1003}{18} - 5\pi^2\Big)\approx 0.05246.
% \end{align*}
\end{enumerate}
\end{theorem}

\section{Relations to population genetics}
\label{6}

Consider sequence data, obtained from a sample of individuals in a
population that reproduces according to Wright-Fisher dynamics.
Besides resampling we consider neutral mutations for an infinite sites
model (as introduced in \cite{Kimura1971}) occurring at rate
$\theta/2$ along each line. Using common notation in population
genetics, we consider the diffusion limit of the dynamics of the
population, where time has been rescaled by a factor $N$, the number
of haploids in the population. The per generation mutation probability
of $\mu$ along each line is rescaled to $\theta/2$, where
$\theta=2N\mu$.

Mutations can also be modelled in the look-down picture: for each
level there is an independent Poisson clock with rate $\theta/2$ by
which mutations on the line carrying the corresponding level
accumulate.  This implies that on each line of the lookdown process
mutations arise at rate $\theta/2$. All mutations an individual
carries at time $t$ are collected along its line of descent.

\subsubsection*{Segregating sites}
For two individuals sampled from the population, the expected number
of segregating sites is $\theta \mathbf E [T_c]$, where $T_c$ is the
random time to coalescence of the individuals' ancestral lineages.
This time is unusually short at instances when the MRCA changes. In
\cite{Tajima1990}, Tajima studied the coalescent at such times. He
concluded that then the coalescence rate from $k$ to $k-1$ ancestral
lineages is $\binom{k+1}{2}$, his argument being that, in addition to
the $k$ lineages, there is one extra line, which apparently must
belong to the family that disappears at the time of the MRCA change.

These coalescence rates can also be seen from the particle
representation of the MRCA process. In fact, the fixation curves give
the shape of the coalescent tree of the whole population back from the
time of the MRCA change.  Recall that the fixation curve moves from
level $k$ to $k+1$ at rate $\binom{k+1}{2}$. Consequently the time
the coalescent back from some time point $E$ stays with $k$ lines is
exponentially distributed with rate $\binom{k+1}{2}$, which means that
this is the rate to go down from $k$ to $k-1$ lineages.  As these
rates differ from the rates in Kingman's coalescent the random
coalescence time $T_c$ of the 2-sample cannot be exponential.
However, the 2-sample coalescent is embedded in the full coalescent;
the probability that the two sampled lines find a common ancestor at
the time when there are $\ell$ lines left in the full coalescent is
(see Remark \ref{remExp} or \cite{Tajima1990}, equation (6))
$$ \frac{2}{(\ell+1)(\ell+2)}.$$ As the time of going down from
infinity to $\ell$ lines in the coalescent at an MRCA time is
distributed like $S_{\ell+1}^\infty$, we obtain the distribution for
the coalescence time $T_c$ of the two lines
$$ \mathbf P[T_c\in dt] = \sum_{\ell=1}^\infty \frac{2}{(\ell+1)(\ell+2)} 
\mathbf P[S_{\ell+1}^\infty\in dt].$$ Taking expectations we obtain
\begin{align*}
 \mathbf E[T_c] &= \sum_{\ell=1}^\infty \frac{2}{(\ell+1)(\ell+2)}
\frac{2}{\ell+1} = 4 \sum_{\ell=1}^\infty \Big( \frac{1}{(\ell+1)^2} -
\frac{1}{(\ell+1)(\ell+2)}\Big) 
%= 4\big( \tfrac{\pi^2}{6} - 1 - \tfrac 12\big) 
= \tfrac 23\pi^2-6 \approx 0.58,
\end{align*}
a result already obtained in \cite{Tajima1990}. As the coalescence
time for two lines in equilibrium is exponential with mean 1, this
result means that the expected number of segregating sites for a
2-sample is reduced by $42\%$ at times when the MRCA changes.

\medskip

For samples of arbitrary size, the number of segregating sites is
Poisson with mean $\theta/2$ times the total branch length of the
sample's genealogical tree. In \cite{RauchBarYam2004}, Figure 2c, a
path of the time evolution of this total branch length is depicted for
a spatial and a ``well-mixed'' population. At certain instances, one
sees sudden substantial decrease of the path length. One may guess
that this happens primarily at times at which the MRCA changes, since
then the coalescent tree is unusually short.

\subsubsection*{Substitutions}
Most mutations that occur in a population are quickly lost. However,
some eventually fixate, i.e. all individuals in the population carry
the new mutation. This replacement is termed a \emph{substitution} and
the corresponding mutations are called \emph{determining mutations}.
In \cite{Watterson1982a} and \cite{Watterson1982b}, Watterson studied
several aspects of the process of substitutions.  While we are
concerned with the jump from today's MRCA to the next one, Watterson
fixes two time points $0$ and $t$ and studies the time between the
MRCAs at these times, i.e. $A_t-A_0$, irrespectively of the number of
MRCAs that are established between $0$ and $t$. All mutations on the
ancestral line between $A_t$ and $A_0$ are then determining mutations
and their number gives the the number of substitutions between times
$0$ and $t$.

The only way a mutation can become a substitution is through an MRCA
change. This is because any mutation that occurs in the population
belongs to one of the two oldest families. For the mutation to become
fixed it is necessary that the family not carrying the mutation dies
out. In other words, it is necessary that the MRCA changes.

Consider the graphical lookdown representation including mutations
falling on lines at all levels at rate $\theta/2$. A mutation that
occurs is determining if and only if it occurs on the line at level
one. Indeed, we already found that MRCAs of the population as seen in
the lookdown picture always are at level one. On the other hand, given
the time point of a mutation on a line at level one, eventually all
individuals in the population are descendants of the individual
carrying this mutation which shows that all mutations that occur at
level one are determining.

Denote by $\mathcal S = \{(\tilde E, S)\}$ the process of times
$\{\tilde E\}$ and number $\{S\}$ of substitutions at these times. As
times $E$ of MRCA changes are the only ones that can be substitution
times and the number of mutations on a line is Poisson distributed
with rate $\tfrac{\theta}{2}$ we find that the process $\mathcal S$ is
a close relative to the MRCA point process:

\begin{proposition}
  Let $\{(E,B)\}$ be distributed as the MRCA point process.
  Additionally, for all successive pairs $(E',B')$ and $(E'',B'')$,
  let $S''$ be Poisson-distributed with intensity parameter
  $\tfrac{\theta}{2}(B'' - B')$.  Then $\{(E'', S''): S''>0 \}$ is a
  version of $\mathcal S$.
\end{proposition}

This confirms the observation in \cite{Watterson1982b} that (i) substitution times do not
form a Poisson process and (ii) substitutions tend to occur in
clusters.

\section{Proof of Theorem \ref{th1}}
\label{7}
Recall the definition of $L_t$ and $I_t$ in \eqref{defL} and
\eqref{defI}, and also recall that we put without loss of generality
$t=0$, omitting the corresponding sub-and superscripts $0$. By
definition, the fixation curve $F_B$ starts at time $B$ at level $2$
and exits at time $E$ at level $\infty$; let us now extend this
definition by putting
\begin{equation*} F_B^\tau := 1 \quad \mbox { if }  \tau < B.
\end{equation*}

 The following auxiliary variables will be helpful:
\begin{equation*}
K^j:= \mbox{ the level of } F_B \mbox{ while }C = j , \qquad j=2,3,\ldots
\end{equation*}
and
\begin{equation*}
I^k:= \mbox{ the level of } C \mbox{ when }F_B \mbox{ reaches level } k, \qquad k=2,3,\ldots
\end{equation*}

\begin{figure}
\begin{center}
\includegraphics[width=13.5cm]{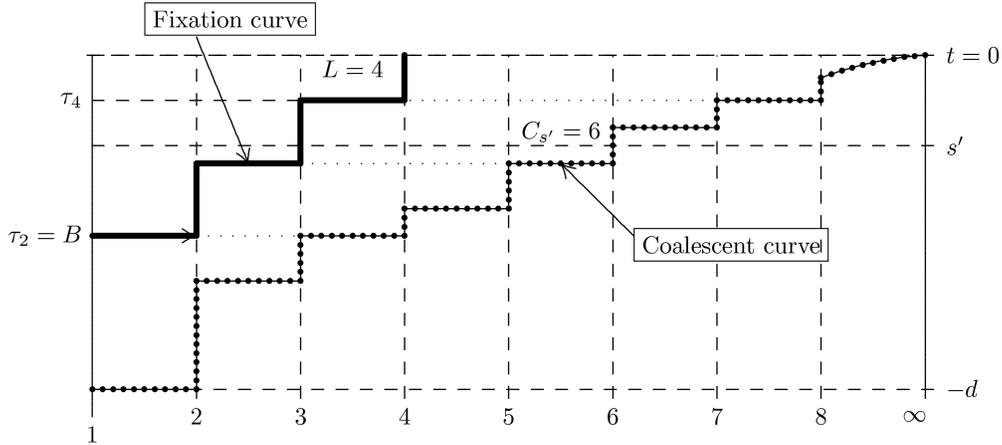}
\end{center}
\caption{\label{FigLookdown2}The variable $I^k$ is the level of the
  coalescence curve back from time $t=0$ when the next fixation curve
  has reached level $k$. The corresponding real times are denoted by
  $\tau_k$. The variable $K^j$ is the level of the fixation curve when
  the coalescent curve has reached level $j$. In this example,
  $I^2=4$, $I^3 = 6$, $I^4=8$, $K^2=K^3=1, K^4=K^5=2, K^6=K^7=3$ and
  $K^8=4$.}
\end{figure}

Formally, putting
$$ \tau_k:=\inf\{\tau: F_B^\tau=k\},\quad k=2,3,..$$ 
we have
$$I^k = \begin{cases} C_{\tau_k}, & \text{ if } \tau_k<0,\\ \infty, &
  \text{ if } \tau_k\geq 0.\end{cases}$$ Thus, $\tau_2 = B$, $I^2 =
I$,
 \begin{equation}\label{LI}
   K^j = \max\{k: I^k \le j\} \quad \text{ and }\quad L = K^\infty:= \lim_{j\to \infty} K^j.
\end{equation}
The random variables $I^k$, $K^j$ and $L$ are illustrated in Figure
\ref{FigLookdown2}.

\begin{lemma}\label{MarkovK}
$K= (K^2, K^3,\ldots)$ is an inhomogeneous Markov chain starting in $K^2=1$ and with transition probability given by
\begin{align}\label{fundamental2}
\mathbf P[K^{j+1}=k+1 | K^j = k] =   \frac{\binom{k+1}{2}}{\binom{j+1}{2}}= 1-\mathbf P[K^{j+1}=k | K^j = k], \quad j > k \ge 1.
\end{align}
Moreover, $K$ is independent of the coalescence curve $C$.
\end{lemma}
\begin{proof}
  When the coalescence curve moves from $j+1$ to $j$ at some time $s$, that is, $C_{s} = j+1$ and $C_{s-} = j$,  then some
  look-down event involving two levels  $\le k+1$ must
happen. When at time $s-$ the  fixation
  curve is at level $k$, the probability that the fixation curve jumps
  at time $s$ from $k$ to $k+1$ is
$\binom{k+1}{2}/\binom{j+1}{2}$,
since all possible look-down events are equally probable and there are
$\binom{\ell+1}{2}$  events that push the next fixation curve one
level up. Observe that this is independent of the prehistory $K^2,\ldots, K^{j-1}$,  independent of the coalescence curve $C$, and in particular also independent of the time $A_0$. 
\end{proof}

In the next lemma we calculate the joint
distribution of the random variables $I^k$. 

\begin{lemma}
\label{L:LD}
The joint distribution of $(I^2,I^3,\ldots)$ is given by
\begin{align}\label{drit}
  \mathbf P[I^2=\ldots=\infty] &= \frac 13,\\ \label{distI}
  \mathbf P[I^2=i_2,\ldots, I^\ell=i_\ell,I^{\ell+1}=\ldots=\infty]&=
  \frac{\ell!(\ell-1)!}{3}\prod_{m=2}^{\ell}
  \frac{1}{(i_m+m)(i_m+m-1)}
\end{align}
for $2<i_2<\ldots < i_\ell.$
\end{lemma}

\begin{proof}
The event $\{I^2=\infty\}$ equals the event that the next fixation curve has not yet started by  time
0, that is the event $\{B>0\} = \{K_j =1 \mbox{ for all } j=1,2,\ldots \}$.  Thus, using \eqref{fundamental2},
\begin{align*}
  \mathbf P[I^2=\infty] & = \prod_{j=3}^\infty \Big(
  1-\frac{1}{\binom{j}{2}}\Big) = \prod_{j=3}^\infty
  \frac{(j+1)(j-2)}{j(j-1)} = \frac{1}{3},
\end{align*}
since the product telescopes. This shows \eqref{drit}. 
To prove \eqref{distI}, we express the event on its left hand side in terms of the variables $K^j$:
\[ \{I^2=i_2,\ldots, I^\ell=i_\ell,I^{\ell+1}=\ldots=\infty\} \] 
\[= \{K^2=\ldots=K^{i_2-1}=1, \ldots, \,K^{i_{\ell-1}}=\ldots=K^{i_\ell-1}=\ell-1,\, K^{i_\ell}=K^{i_\ell+1} = \ldots =\ell\}.
\]
Putting $i_1=2$ and $i_{\ell+1}=\infty$, and using Lemma  \ref{MarkovK} we arrive at
\begin{align*}
  \mathbf P&[I^2=i_2,\ldots, I^\ell=i_\ell, I^{\ell+1}=\infty] =
  \left[ \prod_{m=1}^{\ell} \prod_{j=i_m+1}^{i_{m+1}-1} \Big( 1 -
    \frac{\binom{m+1}{2}}{\binom{j}{2}}\Big) \right]\cdot
  \left[ \prod_{m=1}^{\ell-1} \frac{\binom{m+1}{2}}{\binom{i_{m+1}}{2}} \right] \\
  & = \ell! (\ell-1)! \left[ \prod_{m=1}^{\ell}
    \prod_{j=i_m+1}^{i_{m+1}-1}
    \frac{(j-m-1)(j+m)}{j(j-1)}\right] \left[ \prod_{m=1}^{\ell-1} \frac{1}{(i_{m+1})(i_{m+1}-1)} \right] \\
  & = \ell! (\ell-1)! \left[ \prod_{m=1}^{\ell-1}
    \frac{(i_m-m)\cdots(i_m-1)}{(i_{m+1}-m-1)\cdots
      (i_{m+1}-2)(i_{m+1}-1)}\right.\\& \left.
    \qquad\qquad\qquad\qquad\qquad\qquad\qquad \frac{i_{m+1}\cdots
      (i_{m+1}+m-1)}{(i_m+1)\cdots(i_{m}+m)i_{m+1}}\right]
  \cdot\frac{(i_\ell-\ell)\cdots (i_\ell-1)}{(i_\ell+1)\cdots
    (i_\ell+\ell)} \\&= \frac{\ell! (\ell-1)!}{3}
  \frac{1}{(i_\ell+\ell)(i_\ell+\ell-1)} \left[ \prod_{m=2}^{\ell-1}
    \frac{1}{(i_m+m)(i_m+m-1)} \right] \\& =
  \frac{\ell!(\ell-1)!}{3}\prod_{m=2}^{\ell}
  \frac{1}{(i_m+m)(i_m+m-1)}
\end{align*}
\end{proof}

From the joint distribution of $I^2,I^3,\ldots$ given in Lemma
\ref{L:LD} we obtain because of \eqref{LI} the joint distribution of
$(L,I)$ by projection:

\subsubsection*{Proof of Proposition \ref{LI_prop}}
Because of $\{L=1\} = \{I^2 = \infty\}$, we obtain the assertion of \eqref{IL} for  $\ell =1 $ from \eqref{drit}.
 For $\ell = 2,3,..$ we proceed by induction.
 For $\ell=2$ we have, using again Lemma
\ref{L:LD},
$$\mathbf P[I=i,L=2] = \mathbf P[I^2=i, I^3=\infty] = \frac{2}{3} \frac{1}{(i+2)(i+1)}.$$
If the assertion is true for all $2,\ldots,\ell$, we have
\begin{align*}
  \mathbf P[I=i,\,&L=\ell+1] = \sum_{i<i_3<\ldots < i_{\ell+1}}
  \mathbf P[I^2=i, I^3=i_3,\ldots, I^{\ell+1}=\ell+1,
  I^{\ell+2}=\infty]\\&= \frac{1}{(i+2)(i+1)} \sum_{i<i_3<\ldots <
    i_{\ell+1}} \frac{(\ell+1)!\ell!}{3}
  \prod_{m=3}^{\ell+1} \frac{1}{(i_m+m)(i_m+m-1)}\\
  & = \frac{(\ell+1)\ell}{(i+2)(i+1)} \sum_{i<j}\frac{1}{(j+3)(j+2)}
  \sum_{j+1<i_3<\ldots < i_\ell} \frac{\ell!(\ell-1)!}{3} \prod_{m=3}^\ell \frac{1}{(i_m+m)(i_m+m-1)} \\
  & = \frac{(\ell+1)\ell}{(i+2)(i+1)} \sum_{i<j} \mathbf
  P[I=j+1,L=\ell]
  = \frac{(\ell+1)!\ell(\ell-1)}{3(i+2)(i+1)}\sum_{i<j} \frac{1}{(j+2)\cdots(j+\ell+1)} \\
  & = \frac{(\ell+1)!\ell}{3(i+2)(i+1)}\sum_{i<j} \frac{1}{(j+2)\cdots(j+\ell)} - \frac{1}{(j+3)\cdots(j+\ell+1)}\\
  & = \frac{(\ell+1)!\ell}{3(i+2)(i+1)} \frac{1}{(i+3)\cdots
    (i+\ell+1)}
\end{align*}
and we are done.
\qed

We turn now to the

\subsubsection*{Completion of the Proof of Theorem 1}
Given $\{L=1\} = \{B>0\}$, and independently of $C$ (and therefore also of $A_0$), the time $B$  it takes to enter the next fixation curve is standard exponentially distributed (and therefore distributed like $S_1^2\sim\exp(1)$), and the additional time it takes this fixation curve to exit is distributed like $S_2^\infty$, and is independent of $B$. 

Given $L=\ell \ge 2$, $I=i < \infty$ and  $A_0 = -d$, the time at which the coalescent curve jumps from level $i+1$ to $i$ is distributed like $-R_{i,d}$. By construction, this is also the time $B$ at which the next fixation curve $F_B$ enters. At time $0$, this fixation curve is at level $\ell$; independently of the past, the time it takes until this fixation curve exits is distributed like $S_\ell^\infty$. $\Box$

%The random variable $I$ from Lemma \ref{L:LD} gives the number of
%lines in the coalescent back from time $0$ at the time at which the
%next MRCA lives (recall that $I:= \infty$ if the next MRCA lives after
%time $0$). On the other hand, the random variable $L$ from Lemma
%\ref{L:LD} gives the number of individuals that still have descendants
%at the time when the next MRCA is established. In the case $L=1$ the
%new fixation curve has not started by time $0$, and it takes time
%$S_1^2\sim\exp(1)$ that it enters the population.  Afterwards it takes
%time $S_2^\infty$ for the fixation curve to be completed.

%In the case $L>1$ we have to split the time between $-d$ and 0
%stochastically to obtain $B$. Given $I=i$, this is given by $R_{i,d}$
%as this is the time the coalescent needs to go down to $i$ lines when
%it is known to go down to one line in time $d$.  Additionally when
%$L=\ell$ the new MRCA is found in the population when the line at
%level $\ell+1$ was pushed infinitely often. This occurs after a
%waiting time distributed as $S_\ell^\infty$.

\section{Proof of Theorem \ref{Th:full}}
\label{8}
First we give a formal description of the dynamics of the process
$\Lambda$. Afterwards we derive its equilibrium distribution, and
finally we show that this equilibrium distribution also prevails at
the distinguished times $E$.

\subsubsection{The dynamics of $\Lambda$}
Assume $Z_t=k$ with $L_t^1=\ell_1,\ldots,L_t^k=\ell_k>L_t^{k+1}=1$,
i.e. at time $t$ there are exactly $k$ particles at levels $> 1$; in
other words, exactly $k$ fixation curves are present at time $t$.
Assume level $j$ looks down to level $i$ for $i<j$. If $j>\ell_1+1$
only lines at level greater than $\ell_1+1$ are pushed.  In this case
no particle moves, i.e. $L_t$ stays constant. When $j\leq\ell_1+1$, at
least the level of the next fixation curve increases by one from
$\ell_1$ to $\ell_1+1$ and the corresponding particle moves. The rate
of these events is $\binom{\ell_1+1}{2}$ which equals the rate at
which a fixation curve moves from $\ell_1$ to $\ell_1+1$.
%This is because after this event the individuals at
%the first $\ell_1+1$ levels in the look-down process have a common
%ancestor which is the same as the most recent common ancestor of the
%first $\ell_1$ before the event. 
When $j$ is at most $\ell_2+1$, also the position of the second
fixation curve is increased and the corresponding particle moves.

As look-down events among the first $\ell_1+1$ levels occur at rate
$\binom{\ell_1+1}{2}$, this is also the rate at which the first
particle moves. To be exact, with rate
$\binom{\ell_1+1}{2}-\binom{\ell_2+1}{2}$ only the first particle is
affected, with rate $\binom{l_2+1}{2}-\binom{\ell_3+1}{2}$ the first
two particle move and so on. Additionally, at rate 1, a look-down
event from level $2$ to $1$ occurs which has the effect that a new
particle enters at level 2, i.e. and $L_t^{k+1}$ moves from level $1$
to level $2$ and all particles at levels greater than 1 move as well.

The first particle moves at a quadratic rate and thus reaches infinity
within finite time. When it hits infinity at time $E$ the fixation
curve is completed and $L_E^k = L_{E-}^{k+1}$ for $k\geq 1$ as stated
in \eqref{jumpback} because at time $E$ the second particle becomes
the leading one.

The just stated arguments  prove the following proposition
describing the dynamics of the process $\Lambda$.

\begin{proposition}
\label{propLambda}
 From $\Lambda_t = (\ell_1,\ell_2,\ldots)$,  transitions occur
$$ \begin{cases} \text{ to }(\ell_1+1,\ldots, \ell_k+1, \ell_{k+1},\ldots) 
  \text{ at rate } \binom{\ell_k+1}{2} - \binom{\ell_{k+1}+1}{2}& \text{ if }\ell_k>1 ,\\
  \text{ to }(\ell_1+1,\ldots, \ell_k+1, \ell_{k+1},\ldots) \text{ at
    rate } \qquad \qquad 1 & \text{ if
  }\ell_{k-1}>\ell_k=1.\end{cases}$$
\end{proposition}

To derive the equilibrium distribution of the particle system
$\Lambda$, it will be helpful to compute the one-time distributions
and the limiting distribution of the Markov chain $K$ from Lemma
\ref{MarkovK}.

\begin{lemma}\label{Kdist}
\begin{equation}\label{Kdistformel}
\mathbf P[K^j = k] = \frac{j+1}{j-1} \frac 2{(k+1)(k+2)} , \qquad j > k \ge 1.
\end{equation}
\end{lemma}
\begin{proof}
To prove \eqref{Kdistformel}, we will proceed by induction. Because of $\mathbf P[K^2 = 1]=1$, the formula is true for $j=2$. From \eqref{fundamental2} we obtain the induction step:
\begin{align*}
  \mathbf P[K^{j+1}=k]  &= \mathbf
  P[K^j=k] \Big( 1 -
  \frac{\binom{k+1}{2}}{\binom{j+1}{2}}\Big) + \mathbf
 P[K^j=k-1] \frac{\binom{k}{2}}{\binom{j+1}{2}} \\
 &= \frac{j+1}{j-1} \frac{2}{(k+1)(k+2)}
 \frac{(j+1)j- (k+1)k}{(j+1)j}
 +
  \frac{j+1}{j-1}\frac{2}{k(k+1)}
  \frac{k(k-1)}{(j+1)j} \\
  &= \frac{1}{(j-1)j}
 \frac{2}{(k+1)(k+2)}
 \big((j+1)j - (k+1)k +
(k-1)(k+2)\big) \\ &= \frac{j+2}{j}
 \frac{2}{(k+1)(k+2)}.
\end{align*}
\end{proof}

We are now ready for the

\subsubsection{Completion of the Proof of Theorem 2}
\sloppy We will briefly write $(L^1,L^2,\ldots):=(L^1_0,
L^2_0,\ldots)$. Observe that $L^1$ equals the level $L$ which the
fixation line entering at time $B$ has reached at time $0$. From
\eqref{LI} and Lemma \ref{Kdist} we thus infer readily that
\begin{equation}\label{distL11}
\mathbf P[L_1 = \ell_1] = \frac 2{(\ell_1+1)(\ell_1+2)},
\end{equation}
which proves \eqref{equilL} and also re-establishes \eqref{distL1}.

Next we compute the conditional distribution of $L^{k+1}$, given
$L^k=\ell_k$,\ldots, $L^1=\ell_1$, where $2\le \ell_k < \ldots <
\ell_1$. Consider the $k$-th particle, i.e. the particle which has
level $\ell_k$ at time $0$, and denote the time at which this
particle entered at level 2 by $B_k$.  Since the trajectory of this
particle between times $B_k$ and $0$ is an initial piece of the
coalescent curve (belonging to the exit time of this particle), and
since $L^{k+1}$ is the level of the next fixation curve while this
coalescent curve has level $\ell_k$, we can apply Lemmata
\ref{MarkovK} and \ref{Kdist} to the trajectory of the $(k+1)$-st
particle, parametrised by the levels of the $k$-th particle's
trajectory, to conclude that
\begin{equation}\label{conddist}
\mathbf P[L^{k+1}=\ell_{k+1} | L^k=\ell_k,\ldots, L^1=\ell_1] = 
  \frac{\ell_k+1}{\ell_k-1}\frac{2}{(\ell_{k+1}+1)(\ell_{k+1}+2)}.
  \end{equation}
  Iterating this  we obtain
\begin{equation*}
\begin{aligned}
  \mathbf P&[L^1=\ell_1,\ldots, L^k=\ell_k, L^{k+1}=1] \\& =
  \frac{2}{(\ell_1+1)(\ell_1+2)}
  \frac{\ell_1+1}{\ell_1-1}\frac{2}{(\ell_2+1)(\ell_2+2)}
  \frac{\ell_2+1}{\ell_2-1} \cdots \frac{2}{(\ell_k+1)(\ell_k+2)}
  \frac{\ell_k+1}{\ell_k-1} \frac 13\\& = \frac 13 \prod_{j=1}^k
  \frac{2}{(\ell_j+2)(\ell_j-1)}.
\end{aligned}
\end{equation*}
This shows \eqref {equlambda}.

In Section \ref{7} we argued, by disentangling the combinatorics from the time embedding,  that $L=L_1$ is independent of the coalescence curve $C=C^0$. The same argument shows that $\Lambda_t$ is independent of  $(C^t, \eta_t)$, that is, both the coalescent curve back from time $t$ and the exit time points before $t$. 

We claim that this assertion remains true conditioned on 
      $\{t\in\eta\}$, i.e. the event that $t$ is an exit time. Indeed,
given $\{t \in \eta\}$ we know that $t$ is the exit point of a
fixation curve, which hence must coincide with the coalescence curve
$C^t$. So the above argument  shows that also under this additional
conditioning the particle configuration $\Lambda_t$ is independent of
$C^t$ and $\eta_t$.

\medskip

Now we  turn to assertion 2. of the theorem. Consider a
population in equilibrium.  We know already that $\Lambda_t$ is
in equilibrium, i.e. has distribution $\pi_\Lambda$, independently of
the exit times of particles before $t$ and no matter if $t$ is
conditioned to be an exit time or not.
%As points in $\eta$ after time
%$t$ only depend on $\Lambda_t$, they are independent of $\eta_t$ and
%hence the process $\eta$ is Markov. 
This proves that $\eta$ is a stationary renewal process. Additionally
we know that waiting times between points have the same distribution
as the waiting time out of equilibrium. Thus the waiting times are
memoryless, hence exponential, and $\eta$ is Poisson.  \hfill $\qed$

\begin{remark}\label{quick2}
Here is a more heuristic way (in the spirit of Remark \ref{LIremark}.2) to see  the identity
\begin{equation}\label{L1L2}
\mathbf P[L^{2}=\ell_2 | L^1=\ell_1] = 
  \frac{\ell_1+1}{\ell_1-1}\frac{2}{(\ell_2+1)(\ell_2+2)}.
  \end{equation}
  \end{remark}
Equation \eqref{distL} says that $L=L_1$ has the distribution of the  initial run  length $R$ in a coin tossing with random, uniformly on $[0,1]$ distributed success probability. Similarly, equation \eqref{quick2} says that 
\begin{equation}\label{intuitive}
\mbox{ given } L_1 = \ell_1, \mbox{ the random variable } L_2 \mbox{  is distributed like }R \mbox{ conditioned to }\{R < \ell_1\}.
\end{equation}
This is readily seen because
$$ \mathbf P[R< \ell] = 1 - 2\int_0^1 p^\ell dp = \frac{\ell-1}{\ell+1},$$
and consequently
$$ \mathbf P[R=\ell_2 | R< \ell_1] = \frac{\ell_1+1}{\ell_1-1}\frac{2}{(\ell_2+1)(\ell_2+2)}.$$
The property \eqref{intuitive} can also be understood as follows: 
$L^1=\ell_1$ is the number of currently living
individuals that still have offspring at the time $E_0$ of the next
MRCA change and  $L^2$ is the number of individuals still having offspring at the
time $E_{E_0}$ of the next  but one MRCA change.  At  time $E_0$
one of the two  families which were the oldest at time $0$,   dies out, and our condition is that the $\ell_1$ individuals at time 0 have offspring in the surviving family. 
This surviving family will again be made up of two oldest subfamilies, whose sizes again constitute a uniform split of $[0,1]$. At the time of the next but one MRCA change after time $0$, at least {\em some} of the $\ell_1$ lines must have gone extinct, which amounts to the condition that not all of them belong to the same subfamily. The number $L^2$ of lines that belong to the surviving subfamily thus has the same distribution as $R$ conditioned to $\{R < \ell_1\}$.

\section{Proof of Theorem \ref{LawZ}}
\label{9}
%The independence of $Z$ from the current time of the MRCA follows as
%the look-down picture shows that $Z$ and $Z'$ only depend on the
%distribution of the most persistent lines, 2nd most persistent lines
%and so on in the current coalescent tree and not on coalescent times.

\noindent By the definition of $Z$ we immediately see from Theorem
\ref{Th:full} that in equilibrium
\begin{align}\label{e:Z1}
  \mathbf P[Z=z] & = \sum_{1=\ell_{z+1}<\ell_z<\ldots<\ell_1} \pi_{\Lambda}(\ell_1,\ell_2,\ldots) 
  = \frac 13 \sum_{1<\ell_z<\ldots<\ell_1} \prod_{j=1}^z
  \frac{2}{(\ell_j+2)(\ell_j-1)}
\end{align}
This is the basis for the proof of Theorem \ref{LawZ}. We first show
that the correct weights of the distribution of $Z$ are given by 3.
The weights from 4. are just an application of this. From the weights
we compute the probability generating function given in 1. By
calculating derivatives we obtain the expectation and the variance as
given in 2.

\subsubsection{Proof of 3. and 4.}
All we have to do is to simplify \eqref{e:Z1} for more efficient
computation. Therefore we define
\begin{equation}\label{LawZ3}
\begin{aligned}
f(\ell) &:= \frac{1}{(\ell+2)(\ell-1)},\\
x_k&:= \sum_{\ell=2}^\infty (f(\ell))^k,\\
\end{aligned}
\end{equation}
(the definition of $x_k$ matches the definition in \eqref{Def:xk} as
we will show below) and
\begin{align*}
  p_0&:=1,\\
  p_z & := \sum_{1<\ell_z<\ldots<\ell_1} \prod_{m=1}^z f(\ell_m) =
  \frac{1}{z!}  \sum_{1<\ell_z,\ldots, \ell_1\text{ pwd}} f(\ell_1)
  \cdots f(\ell_{z}).\qquad (z>0)
\end{align*}
Here \emph{pwd} means \emph{pairwise different}. With this definition,
for $z\geq 0$,
\begin{align*}
  \mathbf P[Z=z] = \frac{2^z}{3} p_z.
\end{align*}
We will show first
\begin{align}
  \label{LawZ1}
p_z = \frac 1 z \Big( \sum_{j=1}^z (-1)^{j-1}p_{z-j}x_j \Big),
\end{align}
with $x_j$ given by \eqref{Def:xk}, which gives $p_z$ recursively. Then we we calculate $p_z$ as
\begin{align}
  \label{LawZ2}
  p_z = \sum_{\underline a: \sum i a_i=z} (-1)^{z + \sum a_i}
  \frac{1}{a_1!\cdots a_z!} \prod_{j=1}^z
  \Big(\frac{x_j}{j}\Big)^{a_j},
\end{align}
which gives part 2. of Theorem \ref{LawZ}.

For \eqref{LawZ1} define
\begin{align*}
  b_{z,k} & := \frac{1}{(z-1)!} \sum_{1<\ell_1,\ldots, \ell_z\text{
      pwd}} f(\ell_1) \cdots f(\ell_{z-1})\cdot f(\ell_z)^k.
\end{align*}
Then $p_z =b_{z,1}/z, \quad x_k = b_{1,k}$, \quad
and consequently
\begin{align*}
  b_{z+1,k} &= \frac{1}{z!}\sum_{1<\ell_1,\ldots,\ell_z}
  f(\ell_1)\cdots f(\ell_{z}) \Big( \Big( \sum_{j=2}^\infty
  f(j)^k\Big) - f(\ell_1)^k - \ldots - f(\ell_z)^k\Big) = \frac{1}{z}
  b_{z,1} x_k - b_{z,k+1} \\&= p_z x_k - b_{z,k+1}.
\end{align*}
Therefore we can write
\begin{align*}
  z p_z - (-1)^{z-1} x_z = b_{z,1} + (-1)^{z-2}b_{1,z} =
  \sum_{j=1}^{z-1} (-1)^{j-1} (b_{z+1-j,j} + b_{z-j,j+1}) =
  \sum_{j=1}^{z-1} (-1)^{j-1} p_{z-j} x_j.
\end{align*}
Here the second equality follows because the sum telescopes. This gives 
\eqref{LawZ1}.

The second equation, \eqref{LawZ2} is proved by induction. Instead of
\eqref{LawZ2} we prove
\begin{align}
\label{LawZ2b}
p_z = \sum_{k=1}^z (-1)^{z+k} \frac{1}{k!} \sum_{\underline j:
  j_1+\ldots +j_k=z} \prod_{i=1}^k \frac{x_{j_i}}{j_i}
\end{align}
which then gives \eqref{LawZ2} as the sum is over all vectors
$\underline j$ of length $k$ which sum up to $z$. Every such vector
can be translated into a configuration $\underline a$ with $\sum i
a_i=z$ where $a_i$ is the number of $i$'s in $\underline j$. As for a
given length $k$ of the vector $\underline j$ there are
$\frac{k!}{a_1!\cdots a_k!}$ of these vectors leading to the same
configuration \eqref{LawZ2} is the same as \eqref{LawZ2b}.

For $z=1$ \eqref{LawZ2b} gives $p_1=x_1$ which is true by definition
of $p_z$ and $x_z$. Assume the formula is correct for $1,\ldots,z$ and use
\eqref{LawZ1} to conclude that
\begin{align*}
  p_{z+1} & = \frac{1}{z+1} \Big( \sum_{j=1}^z (-1)^{j-1} x_j
  \sum_{k=1}^{z+1-j} (-1)^{z+1-j+k} \frac{1}{k!}\sum_{\underline
    j:j_1+\ldots + j_k= z+1-j}
  \prod_{i=1}^k \frac{x_{j_i}}{j_i} \Big) + \frac{(-1)^{z}}{z+1} x_{z+1} \\
  & = \Big( \sum_{k=1}^z (-1)^{z+k}\frac{1}{k!} \frac{1}{z+1}
  \sum_{j=1}^{z+1-k} \sum_{\underline j:j_1+\ldots + j_k= z+1-j} x_j
  \prod_{i=1}^k \frac{x_{j_i}}{j_i} \Big) + \frac{(-1)^{z}}{z+1}
  x_{z+1}.
\end{align*}
Since for every $1\leq m\leq k+1$
\begin{align*}
  \sum_{j=1}^{z+1-k} \sum_{\underline j:j_1+\ldots + j_k= z+1-j} x_j
  \prod_{i=1}^k \frac{x_{j_i}}{j_i} = \sum_{\underline j:j_1+\ldots +
    j_{k+1}= z+1} x_{j_m} \prod_{i\neq m} \frac{x_{j_i}}{j_i} =
  \sum_{\underline j:j_1+\ldots + j_{k+1}= z+1} j_m \prod_{i=1}^{k+1}
  \frac{x_{j_i}}{j_i},
\end{align*}
then
\begin{align*}
  \sum_{j=1}^{z+1-k} \sum_{\underline j:j_1+\ldots + j_k= z+1-j} x_j
  \prod_{i=1}^k \frac{x_{j_i}}{j_i} &= \frac{1}{k+1}\sum_{\underline
    j:j_1+\ldots + j_{k+1}= z+1} \prod_{i=1}^{k+1}
  \frac{x_{j_i}}{j_i}\sum_{m=1}^{k+1} j_m \\&= \frac{z+1}{k+1}
  \sum_{\underline j:j_1+\ldots + j_{k+1} = z+1} \prod_{i=1}^{k+1}
  \frac{x_{j_i}}{j_i},
\end{align*}
and therefeore
\begin{align*}
  p_{z+1} & = \Big( \sum_{k=2}^{z+1} (-1)^{z+1+k} \frac{1}{k!}
  \sum_{\underline j:j_1+\ldots + j_{k}= z+1} \prod_{i=1}^k \frac{x_{j_i}}{j_i}\Big) + \frac{(-1)^{z}}{z+1} x_{z+1} \\
  & = \sum_{k=1}^{z+1} (-1)^{z+1+k} \frac{1}{k!}  \sum_{\underline
    j:j_1+\ldots + j_{k}= z+1} \prod_{i=1}^k \frac{x_{j_i}}{j_i},
\end{align*}
which completes the induction and hence proves \eqref{LawZ2b}.

To show that the definition of $x_k$ from \eqref{LawZ3} coincides with
\eqref{Def:xk} we define
\begin{align*}
  A_{k,z} & := \sum_{\ell=2}^\infty \frac{3^{k+z-1}}{(\ell+2)^k
    (\ell-1)^z},
\end{align*}
which gives
\begin{align*}
x_k = \frac{1}{3^{2k-1}}A_{k,k}.
\end{align*}
Thus for $k,z\geq 1, k+z\geq 3$ (otherwise the right side is not
defined)
\begin{align*}
  A_{k,z} = \sum_{\ell=2}^\infty
  \frac{3^{k+z-2}}{(\ell+2)^{k-1}(\ell-1)^{z-1}}\Big( \frac{1}{\ell-1}
  - \frac{1}{\ell+2}\Big) = A_{k-1,z} - A_{k,z-1}.
\end{align*}

Assume we do not sum to $\infty$ but to a large finite $N$ such that
$A_{1,0}$ and $A_{0,1}$ exist. It can be proved by induction on $k+z$
that
\begin{align*}
  A_{k,z} & = (-1)^z \Big(\sum_{j=1}^{k\vee z} \binom{k+z-j-1}{k-j}
  A_{j,0} + (-1)^j \binom{k+z-j-1}{z-j} A_{0,j}\Big) + \mathcal
  O\Big(\frac 1 N\Big)
\end{align*}
where $\binom{-1}{0}=1$. Using this we have
\begin{align*}
  A_{k,k} & = (-1)^k \Big(\sum_{j=1}^{k} \binom{2k-j-1}{k-j}
  \big(A_{j,0} + (1_{\{j\text{ even}\}} - 1_{\{j\text{
      odd}\}})A_{0,j}\big)+ \mathcal O\Big(\frac 1 N\Big).
\end{align*}
So, as
\begin{align*}
  A_{0,j} = 3^{j-1}\zeta(j), \qquad A_{j,0} = 3^{j-1}(\zeta(j) - b_j),
  \qquad b_j = 1+\frac{1}{2^j} + \frac{1}{3^j}
\end{align*}
we can write, now also for $N=\infty$
\begin{align*}
  x_k & = \frac{1}{3^{2k-1}} A_{k,k} = \frac{(-1)^k}{3^{2k-1}} \Big(
  \sum_{j=1}^{k} \binom{2k-j-1}{k-j}3^{j-1}\big(-b_j + 1\{j\text{
    even}\} 2\zeta(j)\big)
\end{align*}
which shows that $x_k$ is of the form \eqref{Def:xk}. This completes the proof of the theorem's assertion 3, from which the weights claimed in assertion 4 follow by inspection.

\subsubsection{Proof of 1.}
To obtain the probability generating function we now calculate
\begin{equation}
\label{sumis1}
\begin{aligned}
  g(t)&:=\mathbf E[t^Z] = \sum_{z=0}^\infty t^z \mathbf P[Z=z] =\frac
  13 \sum_{z=0}^\infty \sum_{k=0}^\infty \sum_{\underline j:
    j_1+\ldots+j_k=z} (2t)^z\frac{(-1)^{z+k}}{k!}\prod_{i=1}^k
  \frac{x_{j_i}}{j_i} \\& = \frac 13 \sum_{k=0}^\infty
  \frac{(-1)^k}{k!} \sum_{j_1,\ldots, j_k} (-2t)^{\sum j_i}
  \prod_{i=1}^k \frac{x_{j_i}}{j_i} \\&= \frac 13 \sum_{k=0}^\infty
  \frac{(-1)^k}{k!}\Big( \sum_{j=1}^\infty \frac{(-2t)^j
    x_j}{j}\Big)^k = \frac 13 \exp\Big( - \sum_{j=1}^\infty
  \frac{(-2t)^j x_j}{j}\Big)
\end{aligned}
\end{equation}
where we have used \eqref{LawZ2b}. The sum in the exponential
simplifies to
\begin{align*}
  -\sum_{j=1}^\infty \frac{(-2t)^j x_j}{j} &= -\sum_{i=2}^\infty
  \sum_{j=1}^\infty \frac 1j \Big( \frac{-2t}{(i+2)(i-1)}\Big)^j =
  \sum_{i=2}^\infty \log\Big( 1 + \frac{2t}{(i+2)(i-1)}\Big)\\&=
  \sum_{i=2}^\infty \log\Big( \frac{i(i+1)+2(t-1)}{(i+2)(i-1)}\Big)
\end{align*}
which proves the formula for the probability generating function.

\subsubsection{Proof of 2.}
We calculate the first two derivatives of the generating function:
\begin{align*}
  g'(t) & = g(t) \sum_{i=2}^\infty\frac{2}{i(i+1)+2(t-1)},\\
  g''(t) & = g(t) \Big(
  \Big(\sum_{i=2}^\infty\frac{2}{i(i+1)+2(t-1)}\Big)^2 -
  \sum_{i=2}^\infty \frac{4}{(i(i+1)+2(t-1))^2}\Big).
\end{align*}
So
\begin{align*}
  \mathbf E[Z] & = g'(1)= 1,\\
  \mathbf{Var}[Z] & = \mathbf E[Z^2]-1 = \mathbf E[Z(Z-1)] = g''(1) =
  1 - 4\sum_{i=2}^\infty \frac{1}{i^2(i+1)^2}
\end{align*}
and the last assertion follows by
\begin{align*}
  \sum_{i=2}^\infty \frac{1}{i^2(i-1)^2} = \sum_{i=2}^\infty \Big(
  \frac{1}{i} - \frac{1}{i+1}\Big)^2 = 2\zeta(2) - (1+1+\tfrac 14) -
  2\tfrac 12 = 2\zeta(2) -\tfrac{13}{4}.
\end{align*}

\subsubsection*{Acknowledgements}
We thank Richard Hudson, John Wakeley and Steve Evans for pointing out
relevant references, and we are  grateful to Tom Kurtz for showing us the recent manuscript \cite{DonnellyKurtz2006} which helped us improve Theorem 2.   Part of our work was done at the Erwin
Schr\"odinger Institute in Vienna, whose hospitality is gratefully
acknowledged.

\bibliographystyle{alpha}

\begin{thebibliography}{STW84}

\bibitem [Bur56]{Burke1956} P.J. Burke. The output of a queueing  
system. {\em Operations Research}, 4 (1956), no. 6, 699-704.

\bibitem[DK96]{DonnellyKurtz1996}
P.~Donnelly and T.G. Kurtz.
\newblock A countable representation of the {F}leming {V}iot measurable
  diffusion.
\newblock {\em Annals of Probability}, 24(2):698--742, 1996.

\bibitem[DK99]{DonnellyKurtz1999}
P.~Donnelly and T.G. Kurtz.
\newblock Particle representations for measure-valued population models.
\newblock {\em Annals of Probability}, 27(1):166--205, 1999.

\bibitem[DK06]{DonnellyKurtz2006}
P.~Donnelly and T.G. Kurtz.
\newblock The Eve Process. Manuscript, personal communication.

\bibitem[Gri80]{Griffiths1980}
R.~C. Griffiths.
\newblock Lines of descent in the diffusion approximation of neutral
  {F}isher-{W}right models.
\newblock {\em Theor. Pop. Biol.}, 17:37--50, 1980.

\bibitem[GT03]{GriffithsTavare2003}
R.~C. Griffiths and S.~Tavar\'e.
\newblock The genealogy of a neutral mutation.
\newblock In {\em Green, P.J., Hjort, N.L. and Richardson, S. (Eds.), Highly
  Structured Stochastic Systems}, pages 393--413. Oxford University Press,
  2003.

\bibitem[JK77]{JohnsonKotz1977}
N.L. Johnson and S.~Kotz.
\newblock {\em Urn {M}odels and their applications}.
\newblock John Wiley \& Sons, 1977.

\bibitem[Kim71]{Kimura1971}
M. Kimura
\newblock Theoretical foundation of population genetics at the molecular level.
\newblock {\em Theo. Pop. Biol.}, 2(2):174--208, 1971.

\bibitem[Kin82]{Kingman1982a}
J.~F.~C. Kingman.
\newblock The coalescent.
\newblock {\em Stochastic Process. Appl.}, 13(3):235--248, 1982.

\bibitem[Kur98]{Kurtz1998}
T.~G. Kurtz.
\newblock Martingale problems for conditional distributions of Markov
  processes.
\newblock {\em Electronic Journal of Probability}, 3, no. 9, pages 1--29, 1998.

\bibitem[Lit75]{Littler1975}
R.A. Littler.
\newblock Loss of variability at one locus in a finite population.
\newblock {\em Math. Biosci.}, 25:151--163, 1975.

%\bibitem[Nor01]{Nordborg2001}
%M. Nordborg.
%\newblock Coalescent Theory. In: Handbook of Statistical Genetics, eds: DJ Balding,. M. Bishop,  C. Cannings, pp 179-212, Wiley  2001.

\bibitem[RBY04]{RauchBarYam2004}
E.M. Rauch and Y. Bar-Yam.
\newblock Theory predicts the uneven distribution of genetic diversity within
  species.
\newblock {\em Nature}, 431:449--452.

\bibitem[STW84]{Saundersetal1984}
I.~W. Saunders, S.~Tavar\'e, and G.~A. Watterson.
\newblock On the genealogy of nested subsamples from a haploid population.
\newblock {\em Adv. Appl. Probab.}, 16:471--491, 1984.

\bibitem[Taj90]{Tajima1990}
F.~Tajima.
\newblock Relationship between {DNA} polymorphism and fixation time.
\newblock {\em Genetics}, 125:447--454, 1990.

\bibitem[Wak05]{Wakeley2005} J.~Wakeley. \newblock Coalescent theory. An Introduction. Roberts \& Company Publishers, Greenwood Village, 2005, and Scion Publishing Ltd (to appear).

\bibitem[Wat82a]{Watterson1982a}
G.A. Watterson.
\newblock Mutant substitutions at linked nucleotide sites.
\newblock {\em Adv. Appl. Prob.}, 14:166--205, 1982.

\bibitem[Wat82b]{Watterson1982b}
G.A. Watterson.
\newblock Substitution times for mutant nucleotides.
\newblock {\em J. Appl. Prob.}, 19A: 59--70, 1982.

\end{thebibliography}

\end{document}